\documentclass[11pt,leqno]{article}
\usepackage{amssymb}
\usepackage[mathscr]{eucal}
\usepackage{amsmath,amssymb,latexsym,theorem,bbm}
\usepackage{color}
\usepackage{appendix,url}
\usepackage{graphicx}

\newcommand{\DD}{\mathbb{D}}

\newcommand{\NN}{\mathbb{N}}

\newcommand{\RR}{\mathbb{R}}

\newcommand{\ZZ}{\mathbb{Z}}

\newcommand{\tA}{\widetilde{A}}

\newcommand{\tC}{\widetilde{C}}

\newcommand{\cA}{{\mathcal A}}

\newcommand{\dd}{\mathrm{d}}
\newcommand{\ee}{\mathrm{e}}

\newcommand{\EE}{\operatorname{\mathbb{E}}}
\newcommand{\PP}{\operatorname{\mathbb{P}}}

\newcommand{\cov}{\operatorname{Cov}}
\newcommand{\sign}{\operatorname{sign}}

\newcommand{\tK}{\widetilde{K}}

\renewcommand{\leq}{\leqslant}
\renewcommand{\geq}{\geqslant}

\newcommand{\proofend}{\hfill\mbox{$\Box$}}

\numberwithin{equation}{section}

\theoremstyle{change} \theorembodyfont{\em}
\newtheorem{Lem}{Lemma.}[section]
\newtheorem{Thm}[Lem]{Theorem.}
\newtheorem{Pro}[Lem]{Proposition.}
\newtheorem{Cor}[Lem]{Corollary.}

\theorembodyfont{\rm}
\newtheorem{Rem}[Lem]{Remark.}

\begin{document}

\begin{center}
 {\bfseries\Large
 Karhunen--Lo\`{e}ve expansion for a generalization of Wiener bridge} \\[5mm]
 {\sc\large
  M\'aty\'as $\text{Barczy}^{*,\diamond}$ and Rezs\H o L. $\text{Lovas}^{**}$}
\end{center}

\vskip0.1cm

\noindent *  MTA-SZTE Analysis and Stochastics Research Group,
             Bolyai Institute, University of Szeged,
             Aradi v\'ertan\'uk tere 1, H--6720 Szeged, Hungary

\noindent ** Institute of Mathematics, University of Debrecen,
              Pf.~400, H--4002 Debrecen, Hungary.

\noindent e--mails: barczy@math.u-szeged.hu (M. Barczy), lovas@science.unideb.hu (R. L. Lovas).

\noindent $\diamond$ Corresponding author.

\renewcommand{\thefootnote}{}
\footnote{ M\'aty\'as Barczy is supported by the J\'anos Bolyai Research
 Scholarship of the Hungarian Academy of Sciences.
Rezs\H o L. Lovas was supported by the Hungarian Scientific Research Fund (OTKA) (Grant No.: K-111651).}
\vspace*{0.2cm}

\vspace*{-10mm}

\begin{abstract}
We derive a Karhunen--Lo\`{e}ve expansion of the Gauss process \ $B_t - g(t)\int_0^1 g'(u)\,\dd B_u$, $t\in[0,1]$,
 \ where \ $(B_t)_{t\in[0,1]}$ \ is a standard Wiener process and \ $g:[0,1]\to \RR$ \ is a twice continuously
 differentiable function with \ $g(0) = 0$ \ and \ $\int_0^1 (g'(u))^2\,\dd u =1$.
\ This process is an important limit process in the theory of goodness-of-fit tests.
We formulate two special cases with the function \ $g(t)=\frac{\sqrt{2}}{\pi}\sin(\pi t)$, $t\in[0,1]$, \
 and \ $g(t)=t$, $t\in[0,1]$, \ respectively.
The latter one corresponds to the Wiener bridge over \ $[0,1]$ \ from \ $0$ \ to \ $0$.
\end{abstract}

\section{Introduction}
In this note we present a new class of Gauss processes, generalizing the Wiener bridge, for which
 Karhunen--Lo\`eve (KL) expansion can be given explicitly.
We point out that there are only few Gauss processes for which the KL expansion is explicitly known.
To give some examples, we mention the Wiener process (see, e.g., Ash and Gardner \cite[Example 1.4.4]{AshGar}),
 the Ornstein--Uhlenbeck process (see, e.g., Papoulis \cite[Problem 12.7]{Papo} or Corlay and Pag\`{e}s \cite[Section 5.4 B]{CorPag}),
 the Wiener bridge (see, e.g., Deheuvels \cite[Remark 1.1]{Deh1}),
 Kac--Kiefer--Wolfowitz process (see, Kac, Kiefer and Wolfowitz \cite{KacKieWol} and Nazarov and Petrova \cite{NazPet}),
 weighted Wiener processes and bridges (Deheuvels and Martynov \cite{DehMar}),
 Jandhyala--MacNeill process (Jandhyala and MacNeill \cite[Section 4]{JanMac}),
 a generalization of Wiener bridge (Pycke \cite{Pyc3}),
 generalized  Anderson--Darling process (Pycke \cite{Pyc1}),
 Rodr\'\i guez--Viollaz process (Pycke \cite{Pyc2}),
 scaled Wiener bridges or also called \ $\alpha$-Wiener bridges (Barczy and Igl\'oi \cite{BarIgl}),
 limit processes related to Cram\'er--von Mises goodness-of-fit tests for hypotheses that
  an observed diffusion process has a sign-type trend coefficient (Gassem \cite{Gas}),
 detrended Wiener processes (Ai, Li and Liu \cite{AiLiLiu}), additive Wiener processes and bridges (Liu \cite{Liu}),
 additive Slepian processes (Liu, Huang and Mao \cite{LiuHuaMao}), Spartan spatial random fields
 (Tsantili and Hristopulos \cite{TsaHri}), the demeaned stationary Ornstein--Uhlenbeck process
 (Ai \cite{Ai}), and the additive two-sided Brownian motion (Ai and Sun \cite{AiSun}).
We also mention that KL expansions of Gauss processes have found several applications
 in small deviation theory, for a complete bibliography, see Lifshits \cite{Lif}.
Here we only mention two papers of Nazarov and Nikitin \cite{Naz1}, \cite{NazNik}.

Let \ $\ZZ_+$, \ $\NN$ \ and \ $\RR$ \ denote the set
 of non-negative integers, positive integers and real numbers, respectively.
For \ $s , t \in \RR$, \ we will use the notation \ $s\land t := \min \{s, t\} $.
\ Let \ $(B_t)_{t\in[0,1]}$ \ be a standard Wiener process,
 and let \ $g:[0,1]\to \RR$ \ be a twice continuously differentiable function such that
 \ $g(0) = 0$ \ and \ $\int_0^1 (g'(u))^2\,\dd u =1$.
Let us introduce the process
 \begin{align}\label{def_g_wiener_bridge}
     Y_t := B_t - g(t)\int_0^1 g'(u)\,\dd B_u,
     \qquad t\in[0,1].
 \end{align}
The process \ $Y=(Y_t)_{t\in[0,1]}$ \ appears as a limit process related to
 a goodness-of-fit test, where one has to decide whether an independent and identically distributed sample
 has a given continuous distribution function depending on some unknown one-dimensional
 parameter, see, e.g., Darling \cite{Dar} or Ben Abdeddaiem \cite{Ben}
 (one can choose \ $h(\theta,t):=g'(t)$, $t\in[0,1]$, \ in formula (1) in \cite{Ben}).
The process \ $Y$ \ also appears as a limit process related to another goodness-of-fit test
 in the case of continuous time observations of a diffusion process with small noise,
 see Ben Abdeddaiem \cite[formula (5)]{Ben}.
One can consider \ $Y$ \ as a generalization of the Wiener bridge corresponding to the function \ $g$.
\ In the special case \ $g(t)=t$, $t\in[0,1]$, \ we have \ $Y_t = B_t - t\int_0^1 1\,\dd B_u = B_t- tB_1$, $t\in[0,1]$,
 \ i.e., it is a Wiener bridge over \ $[0,1]$ \ from \ $0$ \ to \ $0$.
\ However, in general, \ $Y$ \ is not a bridge.
\ Note that \ $Y$ \ is a bridge in the sense that \ $\PP(Y_1 = y_1)=1$ \ with some \ $y_1\in\RR$ \
 (i.e., \ $Y$ \ takes some constant value at time 1 with probability one)
 if and only if \ $g(1)\in\{-1,1\}$, \ and in this case \ $y_1=0$.
\ Indeed, \ $\PP(Y_1 = y_1)=1$ \ with some \ $y_1\in\RR$ \ if and only if
 \ $\DD^2(Y_1)=0$.
\ Since  \ $\DD^2(Y_1) = 1 - (g(1))^2$ \ (see Proposition \ref{Pro1}), we have
 \ $\PP(Y_1 = y_1)=1$ \ with some \ $y_1\in\RR$ \ if and only if
 \ $g(1)\in\{-1,1\}$, \ as desired.
Further, since \ $\EE(Y_1)=0$, \ in this case we have \ $y_1=0$.
\ In the present paper, we do not intend to study whether the process \ $(Y_t)_{t\in[0,1]}$ \ given
 by \eqref{def_g_wiener_bridge} can be considered as a bridge in the sense that it can be derived from
 some appropriate stochastic process (for more information on this procedure, see Barczy and Kern \cite{BarKer}).

We also give a possible {\sl formal} motivation for the definition of the process \ $Y$.
\ Let us write \eqref{def_g_wiener_bridge} in the form
 \[
    \dd Y_t = \dd B_t - \left( \int_0^1 g'(u)\,\dd B_u \right)g'(t)\,\dd t,\qquad t\in[0,1],
 \]
 where \ $\left( \int_0^1 g'(u)\,\dd B_u \right) g'(t)$ \ can be {\sl formally} interpreted as the orthogonal projection
 of the derivative of \ $B_t$ \ (in notation \ $\dd B_t$) onto \ $g'$ \ in \ $L^2$, \ since \ $\int_0^1 (g'(u))^2\,\dd u =1$
 \ (it is only a {\sl formal} one because the derivative of \ $B$ \ does not exist).
So, from this point of view, the derivative of \ $Y_t$ \ (in notation \ $\dd Y_t$) \
 is {\sl formally} the orthogonal component of \ $\dd B_t$ \ with respect to \ $g'$ \ in \ $L^2$,
 \ and one can call \ $\dd Y_t$ \ as the \ $g'$-detrendization of \ $\dd B_t$.

Further, we point out that if \ $g$ \ additionally satisfies \ $g'(1)=0$, \
 then the Gauss process \ $(Y_t)_{t\in[0,1]}$ \ given in \eqref{def_g_wiener_bridge}
 coincides in law with one of the Gauss processes introduced in Nazarov \cite[formula (1.3)]{Naz2},
 for more details, see Appendix \ref{Appendix_Nazarov}.
In the spirit of Nazarov \cite{Naz2}, one can say that \ $(Y_t)_{t\in[0,1]}$ \ is a perturbation of
 the Wiener process \ $(B_t)_{t\in[0,1]}$ \ by the function \ $g$.

\begin{Pro}\label{Pro1}
The process \ $(Y_t)_{t\in[0,1]}$ \ is a zero-mean Gauss process with continuous sample paths almost surely and
 with covariance function \ $R(s,t):=\cov(Y_s,Y_t)  = s\wedge t - g(s)g(t)$, \ $s,t\in[0,1]$.
\end{Pro}

The proof of Proposition \ref{Pro1} can be found in Section \ref{Section_Proofs}.

The continuity of the covariance function \ $R$ \ yields that \ $(Y_t)_{t\in[0,1]}$ \ is \ $L^2$-continuous,
 see, e.g., Theorem 1.3.4 in Ash and Gardner \cite{AshGar}.
We also have \ $R\in L^2([0,1]^2)$.
So, the integral operator associated to the kernel function \ $R$, \ i.e., the operator
 \ $A_R:L^2([0,1])\rightarrow L^2([0,1])$,
\begin{equation}\label{Rdefi}
 \begin{aligned}
  (A_R(\phi))(t):=\int_0^1\! R(t,s)\phi(s)\,\dd s,\quad t\in[0,1],
  \qquad \phi\in L^2[0,1],
\end{aligned}
\end{equation}
 is of the Hilbert--Schmidt type, thus \ $(Y_t)_{t\in[0,1]}$ \ has a Karhunen--Lo\`eve (KL) expansion based on \ $[0,1]$:
\begin{equation} \label{KLrepr}
 Y_t =\sum_{k=1}^\infty
      \sqrt{\lambda_k}\,\xi_k e_k(t),\quad t\in[0,1],
\end{equation}
 where \ $\xi_k,\ k\in\mathbb{N}$, \ are independent standard normally distributed random variables,
 $\lambda_k,\ k\in\NN$, \ are the non-zero (and hence positive) eigenvalues of the integral
 operator \ $A_R$ \ and \ $e_k(t),\ t\in[0,1],\ k\in\NN$, \ are the corresponding normed eigenfunctions,
 which are pairwise orthogonal in \ $L^2([0,1])$, \ see, e.g., Ash and Gardner \cite[Theorem 1.4.1]{AshGar}.
For completeness, we recall that the integral operator \ $A_R$ \ has at most countably many eigenvalues,
 all non-negative (due to positive semi-definiteness) with \ $0$ \ as the only possible limit point, and
 the eigenspaces corresponding to positive eigenvalues are finite dimensional.
Observe that \eqref{KLrepr} has infinitely many terms. Indeed, if it had a finite number of terms,
i.e., if there were only a finite number of eigenfunctions, say $N$, then, by the help of \eqref{def_g_wiener_bridge},
 we would obtain that the Wiener process $(B_t)_{t\in[0,1]}$ is
 concentrated in an at most \ $(N+1)$-dimensional subspace of \ $L^2([0,1])$, \ and so even of \ $C([0,1])$, \ with probability one.
This results in a contradiction, since the integral operator associated to the covariance function
 (as a kernel function) of a standard Wiener process has infinitely many eigenvalues and eigenfunctions.
We also note that the normed eigenfunctions
are unique only up to sign.
The series in \eqref{KLrepr} converges in \ $L^2(\Omega,\cA,\PP)$ \ to \ $Y_t$,
 \ uniformly on \ $[0,1]$, \ i.e.,
 \[
    \sup_{t\in[0,1]}\EE\left( \left\vert Y_t
     - \sum_{k=1}^n\sqrt{\lambda_k}\,\xi_k e_k(t)\right\vert^2\right)\to 0
     \qquad \text{as $n\to\infty$.}
 \]
Moreover, since $R$ is continuous on $[0,1]^2,$  the eigenfunctions corresponding
 to non-zero eigenvalues are also continuous on $[0,1],$ see, e.g. Ash and Gardner \cite[p. 38]{AshGar}
 (this will be important in the proof of Proposition \ref{KLgenthm}, too).
Since the terms on the right-hand side of \eqref{KLrepr} are independent normally distributed random variables
 and $(Y_t)_{t\in[0,1]}$ has continuous sample paths with probability one, the series converges
 even uniformly on $[0,1]$ with probability one (see, e.g., Adler \cite[Theorem 3.8]{Adler}).

\begin{Pro}\label{KLgenthm}
If \ $\lambda$ \ is a non-zero (and hence positive) eigenvalue of the integral operator \ $A_R$ \
 and \ $e$ \ is an eigenfunction corresponding to it, then
 \begin{align}\label{DE_e2}
  \lambda e''(t)= - e(t) - g''(t)\int_0^1 g(s) e(s)\,\dd s,\qquad t\in[0,1],
 \end{align}
 with boundary conditions
 \begin{align}\label{DE_e2_boundary}
   e(0) = 0
     \qquad \text{and}\qquad
    \lambda e'(1) = -g'(1)\int_0^1 g(s)e(s)\,\dd s.
 \end{align}
Conversely, if \ $\lambda$ \ and \ $e(t)$, \ $t\in[0,1]$, \ satisfy \eqref{DE_e2} and \eqref{DE_e2_boundary},
 then \ $\lambda$ \ is an eigenvalue of \ $A_R$ \ and \ $e$ \ is an eigenfunction corresponding to it.
\end{Pro}

Note that for the converse statement in Proposition \ref{KLgenthm} we do not need to know in advance that
 \ $\lambda$ \ is non-zero.
The proof of Proposition \ref{KLgenthm} can be found in Section \ref{Section_Proofs}.

To describe the solutions of \eqref{DE_e2} and \eqref{DE_e2_boundary}, for a fixed $\lambda>0$ we introduce the notations
\begin{gather}\label{help8_constants}
 \begin{split}
a_g(\lambda):=\int_0^1g(t)\cos\left(\frac t{\sqrt\lambda}\right)\dd t,
\qquad b_g(\lambda):=\int_0^1g(t)\sin\left(\frac t{\sqrt\lambda}\right)\dd t,\\
c_g(\lambda):= \int_0^1\left(\int_0^tg(u)g(t)\sin\left(\frac u{\sqrt\lambda}\right)
\cos\left(\frac t{\sqrt\lambda}\right)\dd u \right) \dd t.
\end{split}
\end{gather}

\begin{Thm}\label{THM_KL_gen}
In the KL expansion \eqref{KLrepr} of the process \ $(Y_t)_{t\in[0,1]}$ \ given in \eqref{def_g_wiener_bridge},
 the non-zero (and hence positive) eigenvalues are the solutions of the equation
 \begin{equation}\label{KL_gen_eigenvalue}
 \left(\lambda^{3/2}+\sqrt\lambda\int_0^1g(t)^2\,\dd t+2c_g(\lambda)\right)
 \cos\left(\frac 1{\sqrt\lambda}\right)
 +b_g(\lambda)^2\sin\left(\frac 1{\sqrt\lambda}\right)=0,
 \qquad \lambda>0,
 \end{equation}
 and the corresponding normed eigenfunctions take the form
 \begin{align}\label{KL_gen_eigenvector}
   \begin{split}
  e(t)&=C\Bigg[\sqrt\lambda\cos\left(\frac 1{\sqrt\lambda}\right)g(t)
                +\left(a_g(\lambda)\cos\left(\frac 1{\sqrt\lambda}\right)
                +b_g(\lambda)\sin\left(\frac 1{\sqrt\lambda}\right)\right)
                 \sin\left(\frac t{\sqrt\lambda}\right)\\
      &\phantom{=C\Bigg[\,}  +\cos\left(\frac 1{\sqrt\lambda}\right)
                  \cos\left(\frac t{\sqrt\lambda}\right)
                  \int_0^tg(u)\sin\left(\frac u{\sqrt\lambda}\right)\dd u\\
      &\phantom{=C\Bigg[\,}  -\cos\left(\frac 1{\sqrt\lambda}\right)
                  \sin\left(\frac t{\sqrt\lambda}\right)
                  \int_0^tg(u)\cos\left(\frac u{\sqrt\lambda}\right)\dd u
        \Bigg], \qquad t\in[0,1],
 \end{split}
 \end{align}
 where \ $C\in\RR$ \ is chosen such that \ $\int_0^1 (e(t))^2\,\dd t=1$.
\ (Note that \ $C$ \ may depend on \ $\lambda$, \ but we do not denote this dependence.)
\end{Thm}

The proof of Theorem \ref{THM_KL_gen} can be found in Section \ref{Section_Proofs}.
We emphasize that in Theorem \ref{THM_KL_gen} we give KL expansion \eqref{KLrepr} for a new class of Gauss processes
 with the advantage of an explicit form of the eigenfunctions appearing in \eqref{KLrepr},
 while in the recent papers on KL expansions such as for detrended Wiener processes
 (Ai, Li and Liu \cite{AiLiLiu}), additive Wiener processes and bridges (Liu \cite{Liu})
 and additive Slepian processes (Liu, Huang and Mao \cite{LiuHuaMao}), the form of the eigenfunctions remains somewhat hidden.
As we have already mentioned, in the case of \ $g'(1)=0$, \ the Gauss process \ $(Y_t)_{t\in[0,1]}$ \ coincides in law
 with one of the Gauss processes (1.3) in Nazarov \cite{Naz2}, where he presented a procedure for finding
 the KL expansion for his more general Gauss processes.
In Theorem \ref{THM_KL_gen} we make the KL expansion of \ $(Y_t)_{t\in[0,1]}$ \ as explicit as possible
 by solving the underlying eigenvalue problem directly with the advantage of an explicit form of the eigenfunctions
 unlike in the examples in Section 4 in Nazarov \cite{Naz2}.
We note that Theorem \ref{THM_KL_gen} is applicable in the case of \ $g'(1)\ne 0$ \ as well.

\begin{Rem}\label{Rem_lambda0}
Note that \ $0$ \ may be an eigenvalue of the integral operator \ $A_R$ \
 defined in \eqref{Rdefi}, which is in accordance with Corollary 2 in Nazarov \cite{Naz2}.
For an example, see Section \ref{Section_Proofs}.
\proofend
\end{Rem}

\begin{Rem}
We point out that in the formulation of Theorem \ref{THM_KL_gen} the second derivative
 of \ $g$ \ does not come into play, we use its existence only in the proof of the theorem in question.
This raises the question whether one can find an elementary proof which does not use the existence of \ $g''$, \ only that
 of \ $g'$, \ nor the theory of distributions as in Nazarov \cite[Section 4]{Naz2}.
\ We leave this as an open problem.
The existence of \ $g'$ \ is needed due to the definition of the process \ $Y$, \ see \eqref{def_g_wiener_bridge}.
\proofend
\end{Rem}

In the next remark we recall an application of the KL expansion \eqref{KLrepr}.

\begin{Rem}\label{Rem_Laplace}
The Laplace transform of the \ $L^2([0,1])$-norm square of \ $(Y_t)_{t\in[0,1]}$ \ takes the form
 \begin{align}\label{help_Laplace}
   \EE\left(\exp\left\{-c \int_0^1 Y_t^2\,\dd t \right\}\right)
       = \prod_{k=1}^\infty \frac{1}{\sqrt{1+2c \lambda_k}}\,,
       \qquad c\geq 0.
 \end{align}
Indeed, by \eqref{KLrepr}, we have
 \begin{align*}
    Y_t^2 = \sum_{k=1}^\infty \sum_{\ell=1}^\infty
                           \sqrt{\lambda_k \lambda_\ell}
                           \,\xi_k\xi_\ell\, e_k(t)e_\ell(t),\quad t\in[0,1],
 \end{align*}
 and hence using the fact that \ $(e_k)_{k\in\NN}$ \ is an orthonormal system in \ $L^2([0,1])$, \ we get
 \begin{align*}
   \begin{split}
   \int_0^1  Y_t^2\,\dd t
      = \sum_{k=1}^\infty \sum_{\ell=1}^\infty
        \sqrt{\lambda_k \lambda_\ell} \,\xi_k\xi_\ell
          \int_0^1\!\! e_k(t)e_\ell(t)\,\dd t
      = \sum_{k=1}^\infty \lambda_k\xi_k^2,
  \end{split}
 \end{align*}
 which is nothing else but the Parseval identity in \ $L^2([0,1])$.
\ Since \ $\xi_k$, $k\in\NN$, \ are independent standard normally distributed random variables,
 we get
 \begin{align*}
   \EE\left(\exp\left\{-c \int_0^1 Y_t^2\,\dd t \right\}\right)
         = \prod_{k=1}^\infty \EE\left( \ee^{-c\lambda_k\xi_k^2}\right)
         = \prod_{k=1}^\infty \frac{1}{\sqrt{1+2c\lambda_k}},
         \qquad c\geq 0.
 \end{align*}
\proofend
\end{Rem}

Next we study the special case \ $g(t):=\frac{\sqrt{2}}{\pi}\sin(\pi t)$, $t\in[0,1]$,
 \ yielding \ $\PP(Y_0=0) = \PP(Y_1=B_1)=1$.

\begin{Cor}\label{KLcor1}
If \ $g(t):=\frac{\sqrt{2}}{\pi}\sin(\pi t)$, $t\in[0,1]$, \ then in the KL expansion \eqref{KLrepr} of
 \ $Y_t=B_t - \frac{2}{\pi}\sin(\pi t) \int_0^1 \cos(\pi u)\,\dd B_u$, $t\in[0,1]$, \
 the non-zero (and hence positive) eigenvalues are the solutions of the equation
 \begin{align}\label{eq_eigenvalue}
   \lambda^{3/2} \cos\left(\frac{1}{\sqrt{\lambda}}\right)
      + \frac{2}{\pi^2\left(\pi^2 - \frac{1}{\lambda}\right)}\sin\left(\frac{1}{\sqrt{\lambda}}\right) =0,
      \qquad \lambda\ne \frac{1}{\pi^2},
 \end{align}
 and the corresponding normed eigenfunctions take the form
 \begin{align}\label{eq_eigenfunction}
    e(t) = C \left[\sin\left(\frac{t}{\sqrt{\lambda}}\right)
                    - \frac{2 \sin\left(\frac{1}{\sqrt{\lambda}}\right)}{\lambda\pi\left(\pi^2 - \frac{1}{\lambda}\right)}
                      \sin(\pi t)\right]
         = C \left[\sin\left(\frac{t}{\sqrt{\lambda}}\right)
                    + \sqrt{\lambda}\pi \cos\left(\frac{1}{\sqrt{\lambda}}\right)
                      \sin(\pi t)\right]
 \end{align}
 for \ $t\in[0,1]$, \ where \ $C\in\RR$ \ is chosen such that \ $\int_0^1 (e(t))^2\,\dd t=1$, \ i.e.,
 \[
   C=\pm\left( \frac{\lambda\pi^2}{2} \cos^2\left(\frac{1}{\sqrt{\lambda}}\right)
               + \sqrt{\lambda} \left( \frac{\pi^2}{\pi^2 - \frac{1}{\lambda}} - \frac{1}{4} \right)
                      \sin\left(\frac{2}{\sqrt{\lambda}}\right) + \frac{1}{2}\right)^{-\frac{1}{2}}.
 \]
\end{Cor}

The proof of Corollary \ref{KLcor1} can be found in Section \ref{Section_Proofs}.
In fact, we will provide two proofs.
The first one is an application of Theorem \ref{THM_KL_gen}, which is based on the method of variation of parameters,
 while the second proof is based on the method of undetermined coefficients.

\begin{Rem}\label{Rem_roots}
The equation \eqref{eq_eigenvalue} has a unique root in every interval
 \ $\left(\frac{4}{(2k+1)^2\pi^2}, \frac{4}{(2k-1)^2\pi^2} \right)$, $k\geq 2$, \ $k\in\NN$,
 \ and no root greater than \ $\frac{4}{\pi^2}$.
\ For a proof, see Section \ref{Section_Proofs}.
Since \ $\frac{4}{(2k-1)^2\pi^2}$, $k\in\NN$, \ are the eigenvalues of the integral operator corresponding
 to the covariance function \ $s\wedge t$, $s,t\in[0,1]$, \ of a standard Wiener process,
 we can say that there is a kind of interlacement between the eigenvalues of the integral operators corresponding to the
 underlying standard Wiener process \ $B$ \ and to the perturbed process \ $Y$ \ given in Corollary \ref{KLcor1}.
For more details on this phenomenon in a general setup, see, e.g., Nazarov \cite[page 205]{Naz2}.
Using the rootSolve package in R, we determined the first five roots of  the left-hand side of \eqref{eq_eigenvalue}
 as a function of \ $\lambda>0$, \ listed in decreasing order:
 \[
    \lambda_1=0.338650021, \quad  \frac{1}{\pi^2}\approx  0.101330775,  \quad \lambda_2= 0.021632817,
          \quad \lambda_3= 0.010325434, \quad \lambda_4= 0.006001452.
 \]
The second root we obtained is \ $\frac{1}{\pi^2}$, \ although it is not a solution
 of the equation \eqref{eq_eigenvalue}, neither is it an eigenvalue of \ $A_R$, \ as we shall show
 explicitly in the proof of Corollary 1.7. Its appearance among the roots
 is due to the fact that the left-hand side of \eqref{eq_eigenvalue} may be extended
 continuously to \ $\lambda=\frac{1}{\pi^2}$, \ as we will see in Section 2 in the
 paragraph containing the proofs of the assertions in this remark.
In Figure \ref{fig}, we plotted the left hand side of  \eqref{eq_eigenvalue} as a function of \ $\lambda\in(0,0.35)$.
 \begin{figure}[ht!]
 \centering
 \includegraphics[width=12cm,height=7cm]{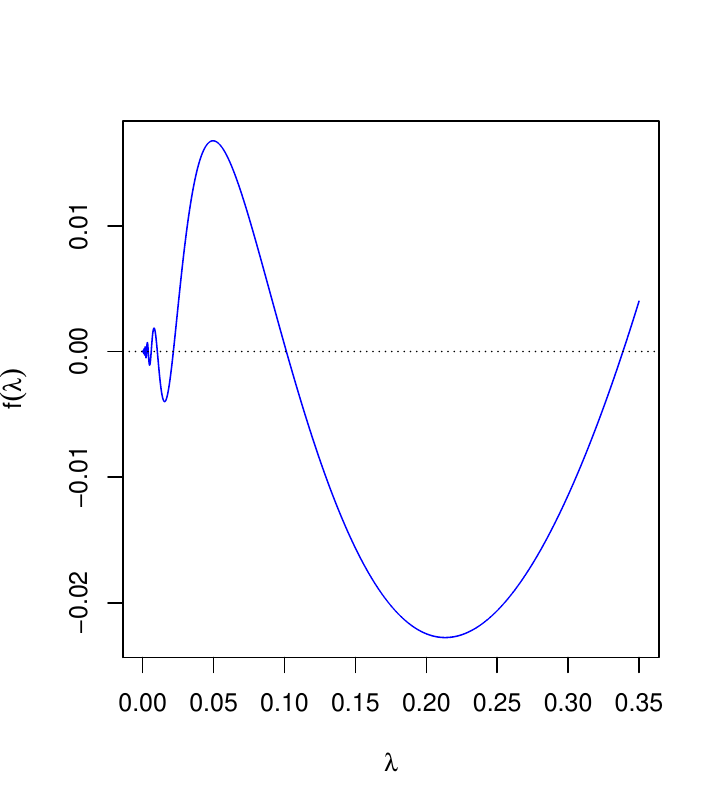}
 \caption{ The function \ $\lambda^{3/2} \cos\left(\frac{1}{\sqrt{\lambda}}\right)
     + \frac{2}{\pi^2(\pi^2 - \frac{1}{\lambda})}\sin\left(\frac{1}{\sqrt{\lambda}}\right)$, \ $\lambda>0$.}
   \label{fig}
 \end{figure}
Hence, by \eqref{help_Laplace}, for small values of \ $c\in(0,1)$, \ the Laplace transform
 \ $\EE\left(\exp\left\{-c \int_0^1 Y_t^2\,\dd t \right\}\right)$ \ of \ $\int_0^1 Y_t^2\,\dd t$, \ can be approximated by
 \begin{align*}
  \Big( 1 + 0.753219c + 0.0530809c^2 + 0.00113549c^3 + 0.000007264c^4 \Big)^{-\frac{1}{2}},
 \end{align*}
 taking into account only the first four terms in the product \eqref{help_Laplace} (corresponding to
 \ $\lambda_1$, \ $\lambda_2$, \ $\lambda_3$ \ and \ $\lambda_4$).
\proofend
\end{Rem}

Finally, we study the special case \ $g(t):=t$, $t\in[0,1]$, \ which is nothing else but the case of a usual Wiener bridge
 over \ $[0,1]$ \ from \ $0$ \ to \ $0$.
\ Note that the KL expansion of a Wiener bridge has been known for a long time, see, e.g., Deheuvels \cite[Remark 1.1]{Deh1}.

\begin{Cor}\label{KLcor2}
If \ $g(t):=t$, $t\in[0,1]$, \ then in the KL expansion \eqref{KLrepr} of the Wiener bridge
 \ $Y_t=B_t - t B_1$, $t\in[0,1]$, \ the non-zero (and hence positive)  eigenvalues are the solutions of the equation
 \begin{align}\label{eq_eigenvalue_WB}
   \sin\left(\frac{1}{\sqrt{\lambda}}\right) =0,
   \qquad \text{i.e.,} \qquad \lambda=\frac{1}{(k\pi)^2}, \quad k\in\NN,
 \end{align}
 and the corresponding normed eigenfunctions take the form
 \begin{align}\label{eq_eigenfunction_WB}
    e(t) = \pm\sqrt{2} \sin(k\pi t), \qquad t\in[0,1],
 \end{align}
 satisfying \ $\int_0^1 (e(t))^2\,\dd t=1$.
\end{Cor}

The proof of Corollary \ref{KLcor2} can be found in Section \ref{Section_Proofs}.

\section{Proofs}\label{Section_Proofs}

\noindent{\bf Proof of Proposition \ref{Pro1}.}
The fact that \ $Y$ \ is a zero-mean Gauss process with continuous sample paths almost surely follows from its definition.
Indeed, since for all \ $0\leq t_1<t_2<\cdots<t_n$, \ $n\in\NN$, \
 \begin{align*}
    \begin{bmatrix}
      Y_{t_1} \\
      \vdots \\
      Y_{t_n} \\
    \end{bmatrix}
   = \begin{bmatrix}
       1 & 0 & \cdots & 0 & -g(t_1) \\
       0 & 1 & \cdots & 0 & -g(t_2) \\
       \vdots & \vdots & \ddots & \vdots & \vdots \\
       0 & 0 & \cdots & 1 & -g(t_n) \\
     \end{bmatrix}
     \begin{bmatrix}
       B_{t_1} \\
       B_{t_2} \\
       \vdots \\
       B_{t_n} \\
       \int_0^1 g'(u)\,\dd B_u \\
     \end{bmatrix},
 \end{align*}
 to check that \ $Y$ \ is a Gauss process it is enough to show that
 \[
    \begin{bmatrix}
      B_{t_1} & \cdots & B_{t_n} & \int_0^1 g'(u)\,\dd B_u \\
    \end{bmatrix}
 \]
 is normally distributed for all \ $0\leq t_1<t_2<\cdots<t_n$, \ $n\in\NN$.
\ This follows from the fact that \ $B$ \ is a Gauss process and from the definition
 of \ $\int_0^1 g'(u)\,\dd B_u$ \ taking into account that an \ $L^2$-limit of normally distributed
 random variables is normally distributed (for a more detailed discussion on a similar procedure,
 see, e.g., the proof of Lemma 48.2 in Bauer \cite{Bau}).
Further, since \ $\int_0^t (g'(u))^2\,\dd u\leq 1$, \ $t\in[0,1]$, \ the process \ $\left(\int_0^t g'(u)\,\dd B_u\right)_{t\in[0,1]}$ \
 is a martingale, and consequently \ $\EE\left(\int_0^1 g'(u)\,\dd B_u\right)=0$, \
 yielding that \ $\EE(Y_t)=0$, \ $t\in[0,1]$.
\ Moreover, for \ $s,t\in[0,1]$,
 \begin{align*}
  R(s,t)
   & = \cov\left(B_s - g(s)\int_0^1 g'(u)\,\dd B_u, B_t - g(t)\int_0^1 g'(u)\,\dd B_u\right) \\
   & = \cov(B_s,B_t)
       - g(t) \cov\left(\int_0^s 1\,\dd B_u, \int_0^1 g'(u)\,\dd B_u\right)\\
   &\phantom{=\;} - g(s) \cov\left(\int_0^1 g'(u)\,\dd B_u, \int_0^t 1\,\dd B_u\right)\\
   &\phantom{=\;}    + g(s)g(t) \cov\left(\int_0^1 g'(u)\,\dd B_u, \int_0^1 g'(u)\,\dd B_u\right)\\
   & = s\wedge t - g(t) \int_0^s g'(u)\,\dd u
       - g(s)\int_0^t g'(u)\,\dd u
       + g(s)g(t) \int_0^1 (g'(u))^2\,\dd u
 \end{align*}
 \begin{align*}
     = s\wedge t - g(t)g(s) - g(s)g(t) + g(s)g(t)
     = s\wedge t - g(s)g(t),
 \end{align*}
 where for the last but one equality we used \ $g(0)=0$ \ and \ $\int_0^1 (g'(u))^2\,\dd u =1$.
\proofend

\bigskip

\noindent{\bf Proof of Proposition \ref{KLgenthm}.}
Let $\lambda$ be a non-zero (and hence positive) eigenvalue of the integral operator \ $A_R$.
\ Then we have
 \begin{equation}\label{eigeq}
 \int_0^1 R(t,s)e(s)\,\dd s=\lambda e(t),\quad t\in[0,1],
 \end{equation}
 and hence
 \begin{align*}
   \int_0^t R(t,s)e(s)\,\dd s + \int_t^1 R(t,s)e(s)\,\dd s
                      =\lambda e(t),\quad t\in[0,1].
 \end{align*}
Then
 \begin{align}\label{help0}
  \begin{split}
  \lambda e(t)
   & = \int_0^t (s-g(s)g(t)) e(s)\,\dd s + \int_t^1 (t-g(s)g(t))e(s)\,\dd s\\
   & = \int_0^t se(s)\,\dd s
      + t\int_t^1 e(s)\,\dd s
      - g(t) \int_0^1 g(s)e(s)\,\dd s,
     \qquad t\in[0,1].
  \end{split}
 \end{align}
The right-hand (and hence the left-hand) side of \eqref{help0} is differentiable with respect to \ $t$,
 \ since \ $e$ \ is continuous (see the Introduction), and, by differentiating with respect to \ $t$, \ we have
 \begin{align*}
   \lambda e'(t)
     &= te(t)
       + \int_t^1 e(s)\,\dd s
       -te(t)
       - g'(t)\int_0^1 g(s)e(s)\,\dd s,
       \qquad t\in[0,1],
 \end{align*}
 yielding that
 \begin{align}\label{DE_e1}
\lambda e'(t)&=-g'(t)\int_0^1 g(s)e(s)\,\dd s
                + \int_t^1 e(s)\,\dd s,\qquad t\in[0,1].
\end{align}
Differentiating \eqref{DE_e1} with respect to \ $t$ \ yields \eqref{DE_e2} (the differentiation is allowed,
 since \ $g$ \ is twice continuously differentiable).
With the special choice \ $t=0$ \ in \eqref{help0}, using that \ $g(0) = 0$ \ and \ $\lambda>0$, \
 we have the boundary condition \ $e(0) = 0$, \ yielding the first part of \eqref{DE_e2_boundary}.
Further, by \eqref{DE_e1} with \ $t=1$, \ we have
 \begin{align*}
   \lambda e'(1) = -g'(1)\int_0^1 g(s)e(s)\,\dd s,
 \end{align*}
 yielding the second part of \eqref{DE_e2_boundary}.

Conversely, let us suppose that \ $\lambda$ \ and \ $e(t)$, $t\in[0,1]$, \ satisfy \eqref{DE_e2} and
 \eqref{DE_e2_boundary}.
Then integration of \eqref{DE_e2} from \ $t$ \ to \ $1$ \ gives
 \[
  \lambda \int_t^1 e''(s)\,\dd s
     = -\int_t^1 e(s)\,\dd s - \int_t^1 g''(s)\,\dd s \int_0^1 g(s)e(s)\,\dd s,
     \qquad t\in[0,1],
 \]
 i.e.,
 \[
   \lambda (e'(1) - e'(t))
     = -\int_t^1 e(s)\,\dd s - (g'(1) - g'(t)) \int_0^1 g(s)e(s)\,\dd s,
          \qquad t\in[0,1].
 \]
By \eqref{DE_e2_boundary}, we have
 \[
   -g'(1)\int_0^1 g(s)e(s)\,\dd s - \lambda e'(t)
       = - \int_t^1 e(s)\,\dd s
         - (g'(1) - g'(t)) \int_0^1 g(s)e(s)\,\dd s,
      \qquad t\in[0,1],
 \]
 i.e.,
 \[
   - \lambda e'(t)
     = - \int_t^1 e(s)\,\dd s + g'(t) \int_0^1 g(s)e(s)\,\dd s,
           \qquad t\in[0,1],
 \]
 which is nothing else but \eqref{DE_e1}.
Integration of \eqref{DE_e1} from \ $0$ \ to \ $t$ \ gives
 \[
   \lambda \int_0^t e'(s)\,\dd s
      = - \int_0^t g'(s)\,\dd s \int_0^1 g(s)e(s)\,\dd s
        + \int_0^t\left(\int_s^1 e(u)\,\dd u\right)\dd s,
        \qquad t\in[0,1],
 \]
 i.e., by integration by parts,
 \[
   \lambda (e(t) - e(0))
      = - (g(t) - g(0)) \int_0^1 g(s)e(s)\,\dd s
        + \int_0^t s e(s)\,\dd s + t\int_t^1 e(s)\,\dd s ,
        \qquad t\in[0,1].
 \]
By \eqref{DE_e2_boundary} and using also \ $g(0)=0$, \ we have
 \begin{align*}
 \lambda e(t)
  & = - g(t) \int_0^1 g(s)e(s)\,\dd s
        + \int_0^t s e(s)\,\dd s + t\int_t^1 e(s)\,\dd s \\
  & = - g(t)\int_0^1 g(s)e(s)\,\dd s
        + \int_0^1 (s\wedge t)e(s)\,\dd s \\
  &  = \int_0^1 ( (s\wedge t) - g(s)g(t) )e(s)\,\dd s
   = \int_0^1 R(t,s) e(s)\,\dd s,
   \qquad t\in[0,1],
 \end{align*}
 i.e., \eqref{eigeq} holds, as desired.
\proofend

\bigskip

\noindent{\bf Proof of Theorem \ref{THM_KL_gen}.}
Let $\lambda>0$ and $e$ be solutions of (1.4) and (1.5), and introduce the notation
\begin{equation}
K:=\int_0^1g(s)e(s)\dd s.\label{eq:Kdef}
\end{equation}
Then (1.4) and (1.5) take the form
\begin{gather}
\lambda e''(t)=-e(t)-Kg''(t),\qquad t\in[0,1],\label{DE_e2withK}\\
e(0)=0\qquad\mbox{and}\qquad\lambda e'(1)=-Kg'(1),\label{eq:boundary}
\end{gather}
respectively.
These are, strictly speaking, not equations for the unknown function $e$ and scalar $\lambda$, since $e$ is hidden also in the coefficient $K$. However, it will prove convenient to consider \eqref{DE_e2withK} temporarily as a second-order linear differential equation (DE) for $e$. The general solution of the homogeneous part of \eqref{DE_e2withK} is
\[
e(t)=c_1\cos\left(\frac t{\sqrt\lambda}\right)
+c_2\sin\left(\frac t{\sqrt\lambda}\right),\qquad t\in[0,1],
\]
where \ $c_1$, $c_2\in\RR$.
To find the solution of the inhomogeneous equation \eqref{DE_e2withK},
 we use the method of variation of parameters, i.e., we are looking for $e$ in the form
\begin{equation}
e(t)=c_1(t)\cos\left(\frac t{\sqrt\lambda}\right)
+c_2(t)\sin\left(\frac t{\sqrt\lambda}\right),\qquad t\in[0,1],\label{eq:varparam}
\end{equation}
with some twice continuously differentiable functions $c_1, c_2:[0,1]\to\RR$.
From this, we obtain the system of equations
\[\begin{split}
\cos\left(\frac t{\sqrt\lambda}\right)c_1'(t)+\sin\left(\frac t{\sqrt\lambda}\right)c_2'(t)&=0\\
      - \sin\left(\frac t{\sqrt\lambda}\right)c_1'(t)
      + \cos\left(\frac t{\sqrt\lambda}\right)c_2'(t)
&=-\frac K {\sqrt{\lambda}}g''(t)
\end{split}\]
for $c_1'(t)$ and $c_2'(t)$. Solving this and substituting the solutions into \eqref{eq:varparam}, we obtain
\begin{align*}
e(t)&=        c_1\cos\left(\frac{t}{\sqrt\lambda}\right)
            + c_2\sin\left(\frac{t}{\sqrt\lambda}\right) \\
   &\quad
     +\frac{K}{\sqrt\lambda}\cos\left(\frac{t}{\sqrt\lambda}\right)
         \int_0^t g''(u)\sin\left(\frac u{\sqrt\lambda}\right)\dd u
     -\frac{K}{\sqrt\lambda}\sin\left(\frac{t}{\sqrt\lambda}\right)
         \int_0^t g''(u)\cos\left(\frac u{\sqrt\lambda}\right)\dd u,
\end{align*}
 where \ $c_1$, $c_2\in\RR$.
 If we take into account the initial condition $e(0)=0$, we can write this in the form
\[
e(t)=c_2\sin\left(\frac{t}{\sqrt\lambda}\right) + \frac{K}{\sqrt\lambda}\cos\left(\frac{t}{\sqrt\lambda}\right)
\int_0^t g''(u)\sin\left(\frac{u}{\sqrt\lambda}\right)\dd u
-\frac{K}{\sqrt\lambda}\sin\left(\frac{t}{\sqrt\lambda}\right)
\int_0^t g''(u)\cos\left(\frac{u}{\sqrt\lambda}\right)\dd u,
\]
where \ $c_2\in\RR$.
\ Applying integration by parts twice in both integrals and taking into account the condition $g(0)=0$, from this we obtain
\[\begin{split}
e(t)&=c_2\sin\left(\frac{t}{\sqrt\lambda}\right)
+\frac{K}{\sqrt\lambda}g'(0)\sin\left(\frac{t}{\sqrt\lambda}\right)-\frac{K}{\lambda}g(t)\\
&\quad+\frac{K}{\lambda^{3/2}}\sin\left(\frac{t}{\sqrt\lambda}\right)
\int_0^t g(u)\cos\left(\frac{u}{\sqrt\lambda}\right)\dd u
-\frac{K}{\lambda^{3/2}}\cos\left(\frac{t}{\sqrt\lambda}\right)
\int_0^t g(u)\sin\left(\frac{u}{\sqrt\lambda}\right)\dd u.
\end{split}\]
With the notation $A:=c_2+  \frac{K}{\sqrt{\lambda}}g'(0)$ the first two terms can be contracted into one:
\begin{equation}\begin{split}
e(t)&=A\sin\left(\frac{t}{\sqrt\lambda}\right)-\frac{K}{\lambda}g(t)
+\frac{K}{\lambda^{3/2}}\sin\left(\frac{t}{\sqrt\lambda}\right)
\int_0^t g(u)\cos\left(\frac{u}{\sqrt\lambda}\right)\dd u\\
&\quad-\frac{K}{\lambda^{3/2}}\cos\left(\frac{t}{\sqrt\lambda}\right)
\int_0^t g(u)\sin\left(\frac{u}{\sqrt\lambda}\right)\dd u.
\end{split}\label{eq:expression_e}\end{equation}
Now we substitute this into the definition \eqref{eq:Kdef} of $K$:
\begin{align*}
 K=\int_0^1g(t)e(t)\dd t
  &=A\int_0^1g(t)\sin\left(\frac{t}{\sqrt\lambda}\right)\dd t-\frac{K}{\lambda}\int_0^1g(t)^2\dd t\\
&\quad +\frac{K}{\lambda^{3/2}}\int_0^1\int_0^t g(u)g(t)
        \cos\left(\frac{u}{\sqrt\lambda}\right)\sin\left(\frac{t}{\sqrt\lambda}\right)\dd u\,\dd t\\
&\quad-\frac{K}{\lambda^{3/2}}\int_0^1\int_0^t g(u)g(t)
       \sin\left(\frac{u}{\sqrt\lambda}\right)\cos\left(\frac{t}{\sqrt\lambda}\right)\dd u\,\dd t,
\end{align*}
where
 \begin{align*}
  \int_0^1&\int_0^t g(u)g(t) \cos\left(\frac{u}{\sqrt\lambda}\right)\sin\left(\frac{t}{\sqrt\lambda}\right)\dd u\,\dd t\\
     &= \int_0^1\int_0^1 g(u)g(t) \cos\left(\frac{u}{\sqrt\lambda}\right)\sin\left(\frac{t}{\sqrt\lambda}\right)\dd u\,\dd t
        - \int_0^1 \left(\int_t^1 g(u)g(t) \cos\left(\frac{u}{\sqrt\lambda}\right)\sin\left(\frac{t}{\sqrt\lambda}\right)\dd u\right)\dd t\\
      &= \int_0^1 g(u) \cos\left(\frac{u}{\sqrt\lambda}\right) \dd u
           \int_0^1 g(t) \sin\left(\frac{t}{\sqrt\lambda}\right) \dd t
      -  \int_0^1 \left(\int_0^u g(u)g(t) \cos\left(\frac{u}{\sqrt\lambda}\right)\sin\left(\frac{t}{\sqrt\lambda}\right)\dd t\right)\dd u\\
       &= a_g(\lambda)b_g(\lambda) - c_g(\lambda).
 \end{align*}
Putting all these together, we obtain the equation
\begin{equation}
b_g(\lambda)A+\left(-1-\frac 1\lambda\int_0^1g(t)^2\dd t
+\frac{a_g(\lambda)b_g(\lambda)-2c_g(\lambda)}{\lambda^{3/2}}\right)K=0.
\label{eq:Kdef1}
\end{equation}
Now we want to substitute $e$ into the second equation of \eqref{eq:boundary}, therefore we calculate the derivative of \eqref{eq:expression_e}:
\begin{align*}
e'(t)&=\frac A{\sqrt\lambda}\cos\left(\frac{t}{\sqrt\lambda}\right)-\frac{K}{\lambda}g'(t)\\
&\quad+\frac{K}{\lambda^2}\cos\left(\frac{t}{\sqrt\lambda}\right)
\int_0^t g(u)\cos\left(\frac{u}{\sqrt\lambda}\right)\dd u
+\frac{K}{\lambda^{3/2}}\sin\left(\frac{t}{\sqrt\lambda}\right)
g(t)\cos\left(\frac{t}{\sqrt\lambda}\right)\\
&\quad+\frac{K}{\lambda^2}\sin\left(\frac{t}{\sqrt\lambda}\right)
\int_0^t g(u)\sin\left(\frac{u}{\sqrt\lambda}\right)\dd u
-\frac{K}{\lambda^{3/2}}\cos\left(\frac{t}{\sqrt\lambda}\right)
g(t)\sin\left(\frac{t}{\sqrt\lambda}\right)\\
&=\frac A{\sqrt\lambda}\cos\left(\frac{t}{\sqrt\lambda}\right)-\frac{K}{\lambda}g'(t)
+\frac{K}{\lambda^2}\cos\left(\frac{t}{\sqrt\lambda}\right)
\int_0^t g(u)\cos\left(\frac{u}{\sqrt\lambda}\right)\dd u\\
&\quad+\frac{K}{\lambda^2}\sin\left(\frac{t}{\sqrt\lambda}\right)
\int_0^t g(u)\sin\left(\frac{u}{\sqrt\lambda}\right)\dd u,\qquad t\in[0,1].
\end{align*}
Substituting this into the second equation of \eqref{eq:boundary}, we get
\[
A\lambda^{3/2}\cos\left(\frac{1}{\sqrt\lambda}\right)
+K\cos\left(\frac{1}{\sqrt\lambda}\right)
\int_0^1 g(u)\cos\left(\frac{u}{\sqrt\lambda}\right)\dd u
+K\sin\left(\frac{1}{\sqrt\lambda}\right)
\int_0^1 g(u)\sin\left(\frac{u}{\sqrt\lambda}\right)\dd u=0,
\]
or, using again the notations \eqref{help8_constants},
\begin{equation}
\lambda^{3/2}\cos\left(\frac{1}{\sqrt\lambda}\right)A
+\left(a_g(\lambda)\cos\left(\frac{1}{\sqrt\lambda}\right)
+b_g(\lambda)\sin\left(\frac{1}{\sqrt\lambda}\right)\right)K=0.\label{eq:boundary1}
\end{equation}
This, together with \eqref{eq:Kdef1}, yields the following homogeneous system of linear equations for the unknowns $A$ and $K$:
\begin{equation}
\begin{split}
b_g(\lambda)A+\left(-1-\frac 1\lambda\int_0^1g(t)^2\dd t
+\frac{a_g(\lambda)b_g(\lambda)-2c_g(\lambda)}{\lambda^{3/2}}\right)K&=0,\\
\lambda^{3/2}\cos\left(\frac{1}{\sqrt\lambda}\right)A
+\left(a_g(\lambda)\cos\left(\frac{1}{\sqrt\lambda}\right)
+b_g(\lambda)\sin\left(\frac{1}{\sqrt\lambda}\right)\right)K&=0.
\end{split}
\label{eq:system}
\end{equation}

In what follows, we show that excluding two special cases, namely, \ $g(t)=t$, $t\in[0,1]$, \ and
 \ $g(t)=-t$, $t\in[0,1]$, \ the function \ $e$ \ given in \eqref{eq:expression_e} can be identically zero
  if and only if \ $A=K=0$.
\ To prove this, it is enough to check that the functions \ $\sin\left(\frac{t}{\sqrt\lambda}\right)$, $t\in[0,1]$, \
 and
 \begin{align}\label{help11}
  \begin{split}
  &-\frac{1}{\lambda}g(t) + \frac{1}{\lambda^{3/2}}\sin\left(\frac{t}{\sqrt\lambda}\right)
                            \int_0^t g(u)\cos\left(\frac{u}{\sqrt\lambda}\right)\dd u \\
  &-\frac{1}{\lambda^{3/2}}\cos\left(\frac{t}{\sqrt\lambda}\right)
          \int_0^t g(u)\sin\left(\frac{u}{\sqrt\lambda}\right)\dd u, \qquad t\in[0,1],
 \end{split}
 \end{align}
 are linearly independent.
On the contrary, let us suppose that they are linearly dependent, i.e., there exist constants \ $\tA,\tK\in\RR$ \
 such that \ $\tA^2+\tK^2\ne 0$ \ and
 \begin{align}\label{help9}
 \begin{split}
  &\tA \sin\left(\frac{t}{\sqrt\lambda}\right)
    +\tK \Bigg( -\frac{1}{\lambda}g(t) + \frac{1}{\lambda^{3/2}}\sin\left(\frac{t}{\sqrt\lambda}\right)
                            \int_0^t g(u)\cos\left(\frac{u}{\sqrt\lambda}\right)\dd u \\
  &\phantom{\tA \sin\left(\frac{t}{\sqrt\lambda}\right) +\tK \Bigg(\,}
    -\frac{1}{\lambda^{3/2}}\cos\left(\frac{t}{\sqrt\lambda}\right)
          \int_0^t g(u)\sin\left(\frac{u}{\sqrt\lambda}\right)\dd u
    \Bigg) = 0 ,\qquad t\in[0,1].
 \end{split}
 \end{align}
By differentiating twice, one can check that
 \begin{align}\label{help10}
   \begin{split}
  \tA \sin\left(\frac{t}{\sqrt\lambda}\right)
     &= -  \tK
        \Bigg( g''(t) - \frac{1}{\lambda} g(t)
                      + \frac{1}{\lambda^{3/2}} \sin\left(\frac{t}{\sqrt\lambda}\right)
                        \int_0^t g(u)\cos\left(\frac{u}{\sqrt\lambda}\right)\dd u\\
     &\phantom{= - \frac{\tK}{\tA} \Bigg(\,}
     - \frac{1}{\lambda^{3/2}} \cos\left(\frac{t}{\sqrt\lambda}\right)
                        \int_0^t g(u)\sin\left(\frac{u}{\sqrt\lambda}\right)\dd u
        \Bigg),
        \qquad t\in[0,1].
    \end{split}
 \end{align}
Comparing \eqref{help9} and \eqref{help10}, we have \ $\tK g''(t) = 0$, $t\in[0,1]$.
\ If \ $\tK=0$, \ then \ $\sin\left(\frac{t}{\sqrt\lambda}\right) = 0$, $t\in[0,1]$, \ which is a contradiction.
Thus \ $\tK\ne 0$, \ and we have \ $g''(t) = 0$, $t\in[0,1]$.
\ Using that \ $g(0) = 0$ \ and \ $\int_0^1(g'(t))^2\,\dd t=1$, \ we get \ $g(t) = t$, $t\in[0,1]$ \
 or \ $g(t) = -t$, $t\in[0,1]$, \ which cases were excluded.
This leads us to a contradiction.

Hence, excluding the two special cases \ $g(t)=t$, $t\in[0,1]$, \ and \ $g(t)=-t$, $t\in[0,1]$, \
 we see that $e$ is not identically zero if and only if at least one of the two coefficients $A$ and $K$ is different from zero, i.e.,
 the system \eqref{eq:system} has a nontrivial solution for $A$ and $K$.
\ This is equivalent to the condition that its determinant is zero, which in turn yields equation \eqref{KL_gen_eigenvalue}.

If \ $g(t) = t$, $t\in[0,1]$, \ or \ $g(t) = -t$, $t\in[0,1]$, \ then, by integration by parts, we have
 \begin{align*}
  & -\frac{1}{\lambda}g(t) + \frac{1}{\lambda^{3/2}}\sin\left(\frac{t}{\sqrt\lambda}\right)
                            \int_0^t g(u)\cos\left(\frac{u}{\sqrt\lambda}\right)\dd u
   -\frac{1}{\lambda^{3/2}}\cos\left(\frac{t}{\sqrt\lambda}\right)
          \int_0^t g(u)\sin\left(\frac{u}{\sqrt\lambda}\right)\dd u \\
  &\quad = \mp \left(\frac{1}{\lambda} \sin\left(\frac{t}{\sqrt\lambda}\right)
               \int_0^t \sin\left(\frac{u}{\sqrt\lambda}\right)\dd u
               + \frac{1}{\lambda} \cos\left(\frac{t}{\sqrt\lambda}\right)
               \int_0^t \cos\left(\frac{u}{\sqrt\lambda}\right)\dd u
         \right) \\
  &\quad = \mp \frac{1}{\sqrt{\lambda}} \sin\left(\frac{t}{\sqrt\lambda}\right),
      \qquad t\in[0,1],
 \end{align*}
 yielding that the function \ $\sin\left(\frac{t}{\sqrt\lambda}\right)$, $t\in[0,1]$, \
 and the function in \eqref{help11} are linearly dependent.
Further, by \eqref{eq:expression_e}, using \ $A=c_2+\frac{K}{\sqrt{\lambda}}g'(0)=c_2\pm\frac{K}{\sqrt{\lambda}}$, \ we have
 \begin{align*}
 &e(t) = \left(A\mp\frac{K}{\sqrt{\lambda}}\right)\sin\left(\frac{t}{\sqrt\lambda}\right)
         = c_2\sin\left(\frac{t}{\sqrt\lambda}\right),\qquad t\in[0,1],\\
 &K=\int_0^1 g(t)e(t)\,\dd t
  = \pm c_2 \int_0^1 t \sin\left(\frac{t}{\sqrt\lambda}\right)\dd t
  = \pm c_2 \sqrt{\lambda} \left( -\cos\left(\frac{1}{\sqrt\lambda}\right)
                                  +\sqrt{\lambda} \sin\left(\frac{1}{\sqrt\lambda}\right) \right) .
 \end{align*}
By the second equation of \eqref{eq:boundary}, we have
 \ $c_2\sqrt{\lambda}\cos\left(\frac{1}{\sqrt\lambda}\right) = \mp K$ \ (which is in fact \eqref{eq:boundary1}
 in the special cases \ $g(t) = t$, $t\in[0,1]$, \ and \ $g(t) = -t$, $t\in[0,1]$).
\ This together with the above form of \ $K$, \ yields \ $\sin\left(\frac{1}{\sqrt{\lambda}}\right)=0$,
 i.e., \ $\lambda = \frac{1}{(k\pi)^2}$, $k\in\NN$.
\ Next we check that the equation \ $\sin\left(\frac{1}{\sqrt{\lambda}}\right)=0$ \ is nothing else but
 the equation \eqref{KL_gen_eigenvalue} in the special cases \ $g(t) = t$, $t\in[0,1]$, \ and
 \ $g(t) = -t$, $t\in[0,1]$.
\ Using integration by parts, the constants defined in \eqref{help8_constants} take the forms
 \begin{align*}
  &a_g(\lambda)=\pm\int_0^1 u\cos\left(\frac{u}{\sqrt{\lambda}}\right)\dd u
               = \pm\left(\sqrt{\lambda}\sin\left(\frac{1}{\sqrt{\lambda}}\right)
                 + \lambda \left( \cos\left(\frac{1}{\sqrt{\lambda}}\right) - 1\right)\right),
                 \qquad \lambda>0,\\
  &b_g(\lambda)=\pm\int_0^1 u\sin\left(\frac{u}{\sqrt{\lambda}}\right)\dd u
               = \pm\left(-\sqrt{\lambda}\cos\left(\frac{1}{\sqrt{\lambda}}\right)
                 + \lambda  \sin\left(\frac{1}{\sqrt{\lambda}}\right)\right),
                 \qquad \lambda>0,\\
  &c_g(\lambda)= \int_0^1 \int_0^t ut  \sin\left(\frac{u}{\sqrt{\lambda}}\right)
                                       \cos\left(\frac{t}{\sqrt{\lambda}}\right) \dd u\,\dd t\\
  &\phantom{c_g(\lambda)\;}
               =-\sqrt{\lambda}\int_0^1 t^2 \cos^2\left(\frac{t}{\sqrt{\lambda}}\right) \dd t
                  +\lambda \int_0^1 t \sin\left(\frac{t}{\sqrt{\lambda}}\right)
                           \cos\left(\frac{t}{\sqrt{\lambda}}\right) \dd t\\
  &\phantom{c_g(\lambda)\;}
               =-\frac{\sqrt{\lambda}}{6} - \frac{\lambda^{3/2}}{2}\cos\left(\frac{2}{\sqrt{\lambda}}\right)
                + \frac{\lambda(\lambda -1)}{4}\sin\left(\frac{2}{\sqrt{\lambda}}\right),\qquad \lambda>0.
 \end{align*}
Hence, using \ $\int_0^1 (g(t))^2\,\dd t = \int_0^1 (\pm t)^2\,\dd t = \frac{1}{3}$, \ the equation \eqref{KL_gen_eigenvalue}
 takes the form
 \begin{align*}
  &\left(\lambda^{3/2} + \frac{\sqrt{\lambda}}{3}
         + 2\left(-\frac{\sqrt{\lambda}}{6} - \frac{\lambda^{3/2}}{2} \cos\left(\frac{2}{\sqrt{\lambda}}\right)
                 + \frac{\lambda(\lambda -1)}{4} \sin\left(\frac{2}{\sqrt{\lambda}}\right) \right)
  \right)
  \cos\left(\frac{1}{\sqrt{\lambda}}\right) \\
  &\qquad+\left( -\sqrt{\lambda}   \cos\left(\frac{1}{\sqrt{\lambda}}\right)
         + \lambda  \sin\left(\frac{1}{\sqrt{\lambda}}\right) \right)^2
    \sin\left(\frac{1}{\sqrt{\lambda}}\right)
    =0.
 \end{align*}
By some algebraic transformations, it is equivalent to \ $\lambda^2 \sin\left(\frac{1}{\sqrt{\lambda}}\right) =0$. \
Since \ $\lambda>0$, \ we have \ $\sin\left(\frac{1}{\sqrt{\lambda}}\right) = 0$, \ yielding
 \ $\lambda=\frac{1}{(k\pi)^2}$, \ $k\in\NN$, \ as desired.

All in all, for every possible \ $g$, \ the equation \eqref{KL_gen_eigenvalue} holds.
It remains to study the form of the eigenfunctions.

If \ $a_g(\lambda)\cos\left(\frac{1}{\sqrt\lambda}\right)+b_g(\lambda)\sin\left(\frac{1}{\sqrt\lambda}\right)\neq 0$,
 \ then from the second equation of \eqref{eq:system} we have
\[
K=-\frac{\lambda^{3/2}\cos\left(\frac{1}{\sqrt\lambda}\right)A}
{a_g(\lambda)\cos\left(\frac{1}{\sqrt\lambda}\right)
+b_g(\lambda)\sin\left(\frac{1}{\sqrt\lambda}\right)}
=-\tC\lambda^{3/2}\cos\left(\frac{1}{\sqrt\lambda}\right),
\]
where
\[
\tC:=\frac A{a_g(\lambda)\cos\left(\frac{1}{\sqrt\lambda}\right)
+b_g(\lambda)\sin\left(\frac{1}{\sqrt\lambda}\right)}.
\]
Finally, if we substitute these expressions for $K$ and $A$ into \eqref{eq:expression_e}, then we obtain \eqref{KL_gen_eigenvector}
 with some appropriately chosen \ $C\in\RR$.

If \ $a_g(\lambda)\cos\left(\frac{1}{\sqrt\lambda}\right)+b_g(\lambda)\sin\left(\frac{1}{\sqrt\lambda}\right)=0$, \ then
 we also show that
  \begin{align}\label{KL_gen_eigenvector_spec}
   \begin{split}
  e(t)&=C\Bigg[\sqrt\lambda\cos\left(\frac 1{\sqrt\lambda}\right)g(t)
                 +\cos\left(\frac 1{\sqrt\lambda}\right)
                  \cos\left(\frac t{\sqrt\lambda}\right)
                  \int_0^tg(u)\sin\left(\frac u{\sqrt\lambda}\right)\dd u\\
      &\phantom{=C\Bigg[\,}  -\cos\left(\frac 1{\sqrt\lambda}\right)
                  \sin\left(\frac t{\sqrt\lambda}\right)
                  \int_0^tg(u)\cos\left(\frac u{\sqrt\lambda}\right)\dd u
        \Bigg], \qquad t\in[0,1],
 \end{split}
 \end{align}
 is a normed eigenvector corresponding to the eigenvalue \ $\lambda$, \
 where \ $C\in\RR$ \ is chosen such that \ $\int_0^1 (e(t))^2\,\dd t=1$.
\ Note that  \eqref{KL_gen_eigenvector_spec} is a special case of \eqref{KL_gen_eigenvector}.
By Proposition \ref{KLgenthm}, it is enough to verify that \ $e(t)$, \ $t\in[0,1]$, \ given
 in \eqref{KL_gen_eigenvector_spec} satisfies \eqref{DE_e2} and \eqref{DE_e2_boundary}.
First, note that, by integration by parts, one can calculate
 \begin{align*}
   \int_0^1 g(s)e(s)\,\dd s
     = C\left[\sqrt{\lambda} \int_0^1 (g(s))^2\,\dd s - a_g(\lambda)b_g(\lambda) + 2c_g(\lambda)\right]
        \cos\left(\frac 1{\sqrt\lambda}\right).
 \end{align*}
Further,
 \begin{align*}
 e'(t) &= C\Bigg[\sqrt{\lambda}\cos\left(\frac 1{\sqrt\lambda}\right) g'(t)
                - \frac{1}{\sqrt{\lambda}}\cos\left(\frac 1{\sqrt\lambda}\right) \sin\left(\frac t{\sqrt\lambda}\right)
                    \int_0^t g(u) \sin\left(\frac u{\sqrt\lambda}\right)\dd u  \\
       &\phantom{= C\Bigg[\;}
                 - \frac{1}{\sqrt{\lambda}}\cos\left(\frac 1{\sqrt\lambda}\right) \cos\left(\frac t{\sqrt\lambda}\right)
                    \int_0^t g(u) \cos\left(\frac u{\sqrt\lambda}\right)\dd u
                    \Bigg], \qquad t\in[0,1],
 \end{align*}
 and
 \begin{align*}
 e''(t) &= C\Bigg[\sqrt{\lambda}\cos\left(\frac 1{\sqrt\lambda}\right) g''(t)
                -  \frac{1}{\sqrt{\lambda}}\cos\left(\frac 1{\sqrt\lambda}\right) g(t)
                - \frac{1}{\lambda}\cos\left(\frac 1{\sqrt\lambda}\right) \cos\left(\frac t{\sqrt\lambda}\right)
                    \int_0^t g(u) \sin\left(\frac u{\sqrt\lambda}\right)\dd u \\
         &\phantom{= C\Bigg[\;}
                + \frac{1}{\lambda}\cos\left(\frac 1{\sqrt\lambda}\right) \sin\left(\frac t{\sqrt\lambda}\right)
                    \int_0^t g(u) \cos\left(\frac u{\sqrt\lambda}\right)\dd u
                    \Bigg], \qquad t\in[0,1],
 \end{align*}
 yielding that
 \begin{align*}
  \lambda e''(t) + e(t) = C\lambda^{3/2} \cos\left(\frac 1{\sqrt\lambda}\right) g''(t),
     \qquad t\in[0,1].
 \end{align*}
Hence, taking into account that \ $C\ne 0$, \ to verify \eqref{DE_e2} it remains to check that
 \[
 \lambda^{3/2} \cos\left(\frac 1{\sqrt\lambda}\right) g''(t)
   = -\left[\sqrt{\lambda} \int_0^1 (g(s))^2\,\dd s - a_g(\lambda)b_g(\lambda) + 2c_g(\lambda)\right]
        \cos\left(\frac 1{\sqrt\lambda}\right)g''(t), \qquad t\in[0,1],
 \]
 which is equivalent to
 \begin{align}\label{help12}
    g''(t) \left[\lambda^{3/2} + \sqrt{\lambda} \int_0^1 (g(s))^2\,\dd s
                  - a_g(\lambda)b_g(\lambda) + 2c_g(\lambda) \right] \cos\left(\frac 1{\sqrt\lambda}\right) = 0,
                  \qquad t\in[0,1].
 \end{align}
Taking into account \eqref{KL_gen_eigenvalue} and
 that \ $a_g(\lambda)\cos\left(\frac{1}{\sqrt\lambda}\right)+b_g(\lambda)\sin\left(\frac{1}{\sqrt\lambda}\right)=0$,
 \ we have
 \begin{align}\label{help13}
     \left[\lambda^{3/2} + \sqrt{\lambda} \int_0^1 (g(s))^2\,\dd s
                  - a_g(\lambda)b_g(\lambda) + 2c_g(\lambda) \right] \cos\left(\frac 1{\sqrt\lambda}\right) = 0,
 \end{align}
 yielding \eqref{help12}.
The boundary conditions \eqref{DE_e2_boundary} hold as well.
Indeed, the boundary condition \ $e(0)=0$ \ is satisfied, since \ $g(0) = 0$, \ and
 the boundary condition \ $\lambda e'(1) = - g'(1)\int_0^1 g(s)e(s)\,\dd s$ \ is equivalent to
 \begin{align*}
  &\cos\left(\frac 1{\sqrt\lambda}\right)
    \Bigg[  g'(1) \left( \lambda^{3/2} +  \sqrt{\lambda} \int_0^1 (g(s))^2\,\dd s
                         - a_g(\lambda)b_g(\lambda) + 2c_g(\lambda)\right) \\
   &\phantom{\cos\left(\frac 1{\sqrt\lambda}\right) \Bigg[}
     - \sqrt{\lambda}\left(a_g(\lambda)\cos\left(\frac{1}{\sqrt\lambda}\right)
       +b_g(\lambda)\sin\left(\frac{1}{\sqrt\lambda}\right)\right) \Bigg] = 0 ,
 \end{align*}
 which is satisfied due to \eqref{help13}.
\proofend

\bigskip

\noindent{\bf An example for the assertion in Remark \ref{Rem_lambda0}.}
Let \ $g:[0,1]\to\RR$, \ $g(t):=\frac{2\sqrt{2}}{\pi}\sin\left(\frac{\pi t}{2}\right)$, $t\in[0,1]$.
\ Then \ $g(0)=0$, \ $g'(1)=0$, \ $\int_0^1(g'(t))^2\,\dd t=1$, \
\ and \ $0$ \ is an eigenvalue of \ $A_R$ \ with \ $-g''(t) = \frac{\pi}{\sqrt{2}}\sin\left( \frac{\pi t}{2}\right)$
 \ as an eigenfunction corresponding to \ $0$, \ which is in accordance with Nazarov \cite[Corollary 2]{Naz2}.
Indeed,
 \begin{align*}
   \int_0^1 R(t,s)e(s)\,\dd s
    & = \int_0^1 \left( s\wedge t - \frac{8}{\pi^2} \sin\left(\frac{\pi s}{2}\right) \sin\left(\frac{\pi t}{2}\right) \right)
          \frac{\pi}{\sqrt{2}}\sin\left(\frac{\pi s}{2}\right) \dd s \\
    & = \frac{\pi}{\sqrt{2}} \Bigg[ \int_0^t s\sin\left(\frac{\pi s}{2}\right)\dd s
                                  + t \int_t^1  \sin\left(\frac{\pi s}{2}\right)\dd s
                                  - \frac{8}{\pi^2} \sin\left(\frac{\pi t}{2}\right)
                                     \int_0^1 \sin^2\left(\frac{\pi s}{2}\right)\dd s
                           \Bigg]
    =0.
 \end{align*}
\proofend

\noindent{\bf First proof of Corollary \ref{KLcor1}.}
We will apply Theorem \ref{THM_KL_gen} with the function \ $g:[0,1]\to\RR$, \ $g(t):=\frac{\sqrt{2}}{\pi}\sin(\pi t)$, \ $t\in[0,1]$.
First, we check that \ $\lambda=\frac{1}{\pi^2}$ \ cannot be an eigenvalue.
On the contrary, let us suppose that \ $\frac{1}{\pi^2}$ \ is an eigenvalue.
Then the constants defined in \eqref{help8_constants} with \ $\lambda=\frac{1}{\pi^2}$ \ take the forms
 \begin{align*}
  &a_g\left( \frac{1}{\pi^2}\right)
     = \frac{\sqrt{2}}{\pi} \int_0^1 \sin(\pi u)\cos(\pi u)\,\dd u =0,
   \qquad
  b_g\left( \frac{1}{\pi^2}\right)
     = \frac{\sqrt{2}}{\pi} \int_0^1 \sin^2(\pi u)\,\dd u =\frac{1}{\sqrt{2}\pi},\\
  &c_g\left( \frac{1}{\pi^2}\right)
     = \frac{2}{\pi^2}\int_0^1\left(\int_0^t \sin^2(\pi u) \sin(\pi t)\cos(\pi t)\,\dd u\right)\dd t
     = -\frac{3}{8\pi^3}.
 \end{align*}
Hence, using \ $\int_0^1 (g(t))^2\,\dd t = \frac{2}{\pi^2}\int_0^1 \sin^2(\pi t)\,\dd t = \frac{1}{\pi^2}$,
 \ \eqref{KL_gen_eigenvalue} with \ $\lambda=\frac{1}{\pi^2}$ \ would imply that
 \begin{align*}
    \left(\frac{1}{\pi^3} + \frac{1}{\pi^3} -\frac{3}{4\pi^3}\right)(-1) + \frac{1}{2\pi^2}\cdot 0=0,
 \end{align*}
 which leads us to a contradiction.

Using \ $\lambda\ne \frac{1}{\pi^2}$ \ and the addition formulas for cosine and sine,
 the constants defined in \eqref{help8_constants} take the forms
 \begin{align*}
  &a_g(\lambda)=\frac{\sqrt{2}}{\pi} \int_0^1 \sin(\pi u)\cos\left(\frac{u}{\sqrt{\lambda}}\right)\dd u
               = \frac{1}{\sqrt{2}\pi}\int_0^1\Bigg[\sin\left( \left(\pi+\frac{1}{\sqrt{\lambda}}\right) u \right)
                                                   + \sin\left( \left(\pi-\frac{1}{\sqrt{\lambda}}\right) u \right)
                                               \Bigg]\dd u \\
  &\phantom{a_g(\lambda)}
              = \frac{\sqrt{2}}{\pi^2 - \frac{1}{\lambda}} \left(1+ \cos\left(\frac{1}{\sqrt{\lambda}}\right)\right),
              \qquad \lambda>0, \quad \lambda\ne\frac{1}{\pi^2},
 \end{align*}
 \begin{align*}
  &b_g(\lambda)=\frac{\sqrt{2}}{\pi} \int_0^1 \sin(\pi u) \sin\left(\frac{u}{\sqrt{\lambda}}\right)\dd u
               = \frac{1}{\sqrt{2}\pi}\int_0^1\Bigg[\cos\left( \left(\pi-\frac{1}{\sqrt{\lambda}}\right) u \right)
                                                    - \cos\left( \left(\pi+\frac{1}{\sqrt{\lambda}}\right) u \right)
                                               \Bigg]\dd u \\
  &\phantom{b_g(\lambda)}
               = \frac{\sqrt{2}}{\pi^2 - \frac{1}{\lambda}} \sin\left(\frac{1}{\sqrt{\lambda}}\right),
                  \qquad \lambda>0, \quad \lambda\ne\frac{1}{\pi^2},
  \end{align*}
  and
  \begin{align*}
  &c_g(\lambda)= \frac{2}{\pi^2} \int_0^1 \int_0^t \sin(\pi u)\sin(\pi t)  \sin\left(\frac{u}{\sqrt{\lambda}}\right)
                                       \cos\left(\frac{t}{\sqrt{\lambda}}\right) \dd u\,\dd t\\
  &\phantom{c_g(\lambda)\;}
               =\frac{1}{\pi^2}
                 \int_0^1 \sin(\pi t) \cos\left(\frac{t}{\sqrt{\lambda}}\right)
                           \left( \frac{ \sin\left( \left(\pi-\frac{1}{\sqrt{\lambda}} \right)t \right)}{\pi-\frac{1}{\sqrt{\lambda}}}
                                  - \frac{ \sin\left( \left(\pi+\frac{1}{\sqrt{\lambda}} \right)t \right)}{\pi+\frac{1}{\sqrt{\lambda}}}
                                 \right)\dd t\\
 &\phantom{c_g(\lambda)\;}
               = \frac{1}{2\pi^2 \left(\pi - \frac{1}{\sqrt{\lambda}}\right)}
                  \int_0^1 \left( \sin\left( \left(\pi+\frac{1}{\sqrt{\lambda}} \right) t \right)
                                  + \sin\left( \left(\pi-\frac{1}{\sqrt{\lambda}} \right) t \right) \right)
                           \sin\left( \left(\pi-\frac{1}{\sqrt{\lambda}} \right) t \right) \dd t \\
 &\phantom{c_g(\lambda)\;=\;}
                - \frac{1}{2\pi^2 \left(\pi + \frac{1}{\sqrt{\lambda}}\right)}
                  \int_0^1 \left( \sin\left( \left(\pi+\frac{1}{\sqrt{\lambda}} \right) t \right)
                                  + \sin\left( \left(\pi-\frac{1}{\sqrt{\lambda}} \right) t \right) \right)
                           \sin\left( \left(\pi+\frac{1}{\sqrt{\lambda}} \right) t \right) \dd t \\
 &\phantom{c_g(\lambda)\;}
               = \frac{1}{2\pi^2}
                \Bigg[ \frac{1}{\pi - \frac{1}{\sqrt{\lambda}}}
                       \int_0^1 \sin^2\left( \left(\pi-\frac{1}{\sqrt{\lambda}} \right) t \right)\dd t
                       - \frac{1}{\pi + \frac{1}{\sqrt{\lambda}}}
                       \int_0^1 \sin^2\left( \left(\pi+\frac{1}{\sqrt{\lambda}} \right) t \right)\dd t \\
 &\phantom{c_g(\lambda)\;= \frac{1}{2\pi^2} \;}
                 + \frac{\frac{2}{\sqrt{\lambda}}}{\pi^2-\frac{1}{\lambda}}
                         \int_0^1 \sin\left( \left(\pi-\frac{1}{\sqrt{\lambda}} \right) t \right)
                                  \sin\left( \left(\pi+\frac{1}{\sqrt{\lambda}} \right) t \right)
                                 \dd t
                \Bigg]\\
 &\phantom{c_g(\lambda)\;}
             = \frac{1}{2\pi^2 \left(\pi^2-\frac{1}{\lambda}\right)}
               \left[\frac{1}{\sqrt{\lambda}} + \frac{\pi^2}{\pi^2-\frac{1}{\lambda}} \sin\left( \frac{2}{\sqrt{\lambda}} \right) \right],
               \qquad \lambda>0, \quad \lambda\ne\frac{1}{\pi^2}.
 \end{align*}
Hence, using \ $\int_0^1 (g(t))^2\,\dd t = \frac{1}{\pi^2}$, \ the equation \eqref{KL_gen_eigenvalue}
 with \ $\lambda\ne \frac{1}{\pi^2}$ \ takes the form
 \begin{align*}
       &\left( \lambda^{3/2}
         + \frac{\sqrt{\lambda}}{\pi^2}
         + \frac{1}{\pi^2 \left(\pi^2 - \frac{1}{\lambda}\right)}
           \left(\frac{1}{\sqrt{\lambda}}  +  \frac{\pi^2}{\pi^2 - \frac{1}{\lambda}}\sin\left( \frac{2}{\sqrt{\lambda}} \right) \right)
           \right)\cos\left( \frac{1}{\sqrt{\lambda}} \right)\\
       &\qquad \qquad + \frac{2}{\left(\pi^2 - \frac{1}{\lambda}\right)^2} \sin^3\left( \frac{1}{\sqrt{\lambda}} \right) =0,
       \qquad \lambda>0, \quad \lambda\ne \frac{1}{\pi^2}.
 \end{align*}
  and, by some algebraic transformations, we have
 \[
  \frac{\lambda^{3/2}}{\pi^2 - \frac{1}{\lambda}} \cos\left( \frac{1}{\sqrt{\lambda}} \right)
     + \frac{2}{\pi^2\left( \pi^2 - \frac{1}{\lambda}\right)^2} \sin\left( \frac{1}{\sqrt{\lambda}} \right)
     = 0, \qquad \lambda>0, \;\; \lambda\ne \frac{1}{\pi^2}.
 \]
By multiplying this equation with \ $\pi^2 - \frac{1}{\lambda}$, \ we obtain that the equation \eqref{KL_gen_eigenvalue}
 with \ $\lambda\ne \frac{1}{\pi^2}$ \ is equivalent to \eqref{eq_eigenvalue}.

Further, the normed eigenfunctions \eqref{KL_gen_eigenvector} with \ $\lambda\ne \frac{1}{\pi^2}$ \ take the form
 \begin{align*}
  e(t) & = C\Bigg[ \frac{\sqrt{2}}{\pi} \sqrt{\lambda} \cos\left( \frac{1}{\sqrt{\lambda}} \right) \sin(\pi t)
                 + \frac{\sqrt{2}}{\pi^2 - \frac{1}{\lambda}}
                   \left( \left( 1 +  \cos\left( \frac{1}{\sqrt{\lambda}} \right)\right) \cos\left( \frac{1}{\sqrt{\lambda}} \right)
                                + \sin^2\left( \frac{1}{\sqrt{\lambda}} \right)  \right)
                  \sin\left( \frac{t}{\sqrt{\lambda}} \right)\\
  & \phantom{= C\Bigg[}
       + \frac{\sqrt{2}}{\pi}\cos\left( \frac{1}{\sqrt{\lambda}} \right) \cos\left( \frac{t}{\sqrt{\lambda}} \right)
                    \int_0^t \sin(\pi u)\sin\left( \frac{u}{\sqrt{\lambda}} \right)\dd u\\
  & \phantom{= C\Bigg[}
      - \frac{\sqrt{2}}{\pi}\cos\left( \frac{1}{\sqrt{\lambda}} \right) \sin\left( \frac{t}{\sqrt{\lambda}} \right)
                    \int_0^t \sin(\pi u)\cos\left( \frac{u}{\sqrt{\lambda}} \right)\dd u
          \Bigg]\\
 &= C\Bigg[ \frac{\sqrt{2}}{\pi} \sqrt{\lambda} \cos\left( \frac{1}{\sqrt{\lambda}} \right) \sin(\pi t)
                 + \frac{\sqrt{2}}{\pi^2 - \frac{1}{\lambda}}
                   \left( 1+ \cos\left( \frac{1}{\sqrt{\lambda}} \right)\right)
                    \sin\left( \frac{t}{\sqrt{\lambda}} \right) \\
  & \phantom{= C\Bigg[}
       + \frac{1}{\sqrt{2}\pi}\cos\left( \frac{1}{\sqrt{\lambda}} \right) \cos\left( \frac{t}{\sqrt{\lambda}} \right)
          \left( \frac{ \sin\left( \left( \pi - \frac{1}{\sqrt{\lambda}}\right)t\right) }{\pi-\frac{1}{\sqrt{\lambda}}}
                 - \frac{ \sin\left( \left( \pi + \frac{1}{\sqrt{\lambda}}\right)t\right) }{\pi+\frac{1}{\sqrt{\lambda}}}
          \right)\\
  & \phantom{= C\Bigg[}
      - \frac{1}{\sqrt{2}\pi}\cos\left( \frac{1}{\sqrt{\lambda}} \right) \sin\left( \frac{t}{\sqrt{\lambda}} \right)
        \left( - \frac{ \cos\left( \left( \pi + \frac{1}{\sqrt{\lambda}}\right)t\right) }{\pi+\frac{1}{\sqrt{\lambda}}}
               - \frac{ \cos\left( \left( \pi - \frac{1}{\sqrt{\lambda}}\right)t\right) }{\pi-\frac{1}{\sqrt{\lambda}}}
               + \frac{1}{\pi+\frac{1}{\sqrt{\lambda}}} + \frac{1}{\pi-\frac{1}{\sqrt{\lambda}}}
          \right) \Bigg]\\
   &= C\Bigg[ \frac{\sqrt{2}}{\pi} \sqrt{\lambda} \cos\left( \frac{1}{\sqrt{\lambda}} \right) \sin(\pi t)
              + \frac{\sqrt{2}}{\pi^2 - \frac{1}{\lambda}} \sin\left( \frac{t}{\sqrt{\lambda}} \right)\\
    & \phantom{= C\Bigg[}
              + \frac{1}{\sqrt{2}\pi \left(\pi - \frac{1}{\sqrt{\lambda}}\right)}
                \cos\left( \frac{1}{\sqrt{\lambda}} \right)
                \left( \sin\left( \left( \pi - \frac{1}{\sqrt{\lambda}}\right)t\right)
                       \cos\left( \frac{t}{\sqrt{\lambda}} \right)
                       + \cos\left( \left( \pi - \frac{1}{\sqrt{\lambda}}\right)t\right)
                         \sin\left( \frac{t}{\sqrt{\lambda}} \right)
                 \right)\\
  & \phantom{= C\Bigg[}
              + \frac{1}{\sqrt{2}\pi \left(\pi + \frac{1}{\sqrt{\lambda}}\right)}
                \cos\left( \frac{1}{\sqrt{\lambda}} \right)
                \left( \sin\left( \frac{t}{\sqrt{\lambda}} \right)
                        \cos\left( \left( \pi + \frac{1}{\sqrt{\lambda}}\right)t\right)
                       - \cos\left( \frac{t}{\sqrt{\lambda}} \right)
                         \sin\left( \left( \pi + \frac{1}{\sqrt{\lambda}}\right)t\right)
                 \right)\Bigg]\\
  & = C\Bigg[ \frac{\sqrt{2}}{\pi}\sqrt{\lambda} \cos\left( \frac{1}{\sqrt{\lambda}} \right) \sin(\pi t)
              + \frac{\sqrt{2}}{\pi^2 - \frac{1}{\lambda}} \sin\left( \frac{t}{\sqrt{\lambda}} \right)\\
  & \phantom{= C\Bigg[}
              + \frac{1}{\sqrt{2}\pi \left(\pi - \frac{1}{\sqrt{\lambda}}\right)}
                 \cos\left( \frac{1}{\sqrt{\lambda}} \right)
                 \sin(\pi t)
               + \frac{1}{\sqrt{2}\pi \left(\pi + \frac{1}{\sqrt{\lambda}}\right)}
                 \cos\left( \frac{1}{\sqrt{\lambda}} \right)
                 \sin(-\pi t)
       \Bigg] \\
  & = \frac{C\sqrt{2}}{\pi^2 - \frac{1}{\lambda}}
      \left( \sin\left( \frac{t}{\sqrt{\lambda}} \right)
             +\sqrt{\lambda}\pi \cos\left( \frac{1}{\sqrt{\lambda}} \right) \sin(\pi t) \right),
 \end{align*}
 where \ $C\in\RR$ \ is such that \ $\int_0^1 (e(t))^2\,\dd t=1$.
\ Hence
 \begin{align*}
  \frac{1}{C^2}
   & = \frac{2}{\left(\pi^2 - \frac{1}{\lambda}\right)^2}
      \int_0^1 \left( \sin^2\left( \frac{t}{\sqrt{\lambda}} \right)
                      + \lambda\pi^2 \cos^2\left( \frac{1}{\sqrt{\lambda}} \right) \sin^2(\pi t)
                      + 2\sqrt{\lambda}\pi \cos\left( \frac{1}{\sqrt{\lambda}} \right)
                          \sin\left( \frac{t}{\sqrt{\lambda}} \right)
                          \sin(\pi t)
                \right)\dd t\\
   & = \frac{2}{\left(\pi^2 - \frac{1}{\lambda}\right)^2}
       \left( \frac{\lambda\pi^2}{2} \cos^2\left(\frac{1}{\sqrt{\lambda}}\right)
                + \sqrt{\lambda} \left( \frac{\pi^2}{\pi^2 - \frac{1}{\lambda}} - \frac{1}{4} \right)
                   \sin\left(\frac{2}{\sqrt{\lambda}}\right) + \frac{1}{2} \right) .
 \end{align*}
Merging the factor \ $\frac{\sqrt{2}}{\pi^2 - 1/\lambda}$ \ into \ $C$,
 \ taking into account \eqref{eq_eigenvalue} and that
 \[
    \sqrt{\lambda} \pi \cos\left(\frac{1}{\sqrt{\lambda}}\right)
       = -\frac{2}{\lambda \pi \left(\pi^2 - \frac{1}{\lambda}\right) }
          \sin\left(\frac{1}{\sqrt{\lambda}}\right),
          \qquad  \lambda>0, \quad \lambda\ne \frac{1}{\pi^2},
 \]
 this yields \eqref{eq_eigenfunction}.
\proofend

\bigskip

\noindent{\bf Second proof of Corollary \ref{KLcor1}.}
In the special case \ $g(t):=\frac{\sqrt{2}}{\pi}\sin(\pi t)$, $t\in[0,1]$, \
 the DE \eqref{DE_e2} and the boundary conditions \eqref{DE_e2_boundary} take the form
 \begin{align}\label{DE_e2_spec}
    \lambda e''(t)
       = -e(t) + 2\sin(\pi t) \int_0^1 \sin(\pi s) e(s)\,\dd s,
       \qquad t\in[0,1],
 \end{align}
 and
 \[
   e(0)=0 \qquad \text{and}\qquad \lambda e'(1) = \frac{2}{\pi} \int_0^1 \sin(\pi s)e(s)\,\dd s,
 \]
 respectively.
With the special choice \ $t=0$, \ using \ $e(0) = 0$ \ and \ $\lambda>0$,
 \ we have \ $e''(0) = 0$.
The DE \eqref{DE_e2_spec} is a second order linear inhomogeneous DE
 of the type \ $\lambda e''(t) = - e(t) + B\sin(\pi t)$, $t\in[0,1]$,
 \ where \ $B:=2\int_0^1 \sin(\pi s)e(s)\,\dd s$ \ (for more details, see the proof of Theorem \ref{THM_KL_gen}).
\ By the method of undetermined coefficients, its general solution takes the form
 \begin{align}\label{help2}
   e(t) = a\sin\left(\frac{t}{\sqrt{\lambda}}\right)
          + b\cos\left(\frac{t}{\sqrt{\lambda}}\right)
          + c\sin(\pi t)
          + d\cos(\pi t),
         \qquad t\in[0,1],
 \end{align}
 where \ $a,b,c,d\in\RR$.
\ Hence
 \begin{align}\label{help3}
   e'(t) & = \frac{a}{\sqrt{\lambda}}\cos\left(\frac{t}{\sqrt{\lambda}}\right)
              - \frac{b}{\sqrt{\lambda}} \sin\left(\frac{t}{\sqrt{\lambda}}\right)
              + c\pi\cos(\pi t)
              - d\pi\sin(\pi t),
            \qquad t\in[0,1],\\ \label{help4}
   e''(t) & = - \frac{a}{\lambda}\sin\left(\frac{t}{\sqrt{\lambda}}\right)
              - \frac{b}{\lambda} \cos\left(\frac{t}{\sqrt{\lambda}}\right)
              - c\pi^2\sin(\pi t)
              - d\pi^2\cos(\pi t),
            \qquad t\in[0,1].
 \end{align}
Then the DE \eqref{DE_e2_spec} takes the form
 \begin{align}\label{help5}
   c(1-\lambda \pi^2) \sin(\pi t)  + d(1-\lambda \pi^2)\cos(\pi t)
      = 2\sin(\pi t) \int_0^1 \sin(\pi s) e(s)\,\dd s,
      \qquad t\in[0,1].
 \end{align}
Since \ $e(0) = 0$, \ by \eqref{help2}, we have \ $d=-b$.
Since \ $e''(0) = 0$, \ by \eqref{help4}, we have \ $0 = -\frac{b}{\lambda} - d\pi^2$, \ and then,
 since \ $d=-b$, \ we get \ $b=d=0$ \ or \ $\lambda=\frac{1}{\pi^2}$.

Next, we check that \ $\lambda=\frac{1}{\pi^2}$ \ cannot be an eigenvalue.
On the contrary, let us suppose that \ $\lambda=\frac{1}{\pi^2}$ \ is an eigenvalue.
Then, since \ $d=-b$, \ by \eqref{help2}, we have
 \[
  e(t) = a \sin(\pi t) + c\sin(\pi t)
       = (a+c)\sin(\pi t) ,\qquad t\in[0,1],
 \]
 where \ $a,c\in\RR$.
\ Further, by \eqref{help5},
 \[
  0 = \int_0^1 \sin(\pi s) e(s)\,\dd s
    = (a+c) \int_0^1 (\sin(\pi s))^2 \,\dd s
    = \frac{a+c}{2},
 \]
 yielding \ $e(t) = 0$, $t\in[0,1]$, \ which leads us to a contradiction.

Since \ $\lambda\ne\frac{1}{\pi^2}$, \ then we have \ $b=d=0$.
\ By \eqref{help2} and \eqref{help0} together with \ $g(1) = \frac{\sqrt{2}}{\pi}\sin(\pi)=0$, \ we have
 \ $e(t)=a\sin\left(\frac{t}{\sqrt{\lambda}}\right) + c\sin(\pi t)$, \ $t\in[0,1]$, \
 and \ $\lambda e(1) = \int_0^1 se(s)\,\dd s$.
\ Hence
 \begin{align*}
  \lambda a\sin\left(\frac{1}{\sqrt{\lambda}}\right)
   & = a\int_0^1 s \sin\left(\frac{s}{\sqrt{\lambda}}\right)\dd s
      + c \int_0^1 s \sin(\pi s)\dd s  \\
   & = a\left( -\sqrt{\lambda}\cos\left(\frac{1}{\sqrt{\lambda}}\right)
               + \lambda \sin\left(\frac{1}{\sqrt{\lambda}}\right) \right)
      +\frac{c}{\pi}.
 \end{align*}
Then
 \begin{align}\label{help6}
   -a\sqrt{\lambda} \cos\left(\frac{1}{\sqrt{\lambda}}\right) + \frac{c}{\pi} = 0,
   \qquad \lambda>0, \quad \lambda\ne \frac{1}{\pi^2}.
 \end{align}
By \eqref{help3} and \ $\lambda e'(1) = \frac{2}{\pi}\int_0^1\sin(\pi s)e(s)\,\dd s$, \ we have
 \begin{align*}
 a\sqrt{\lambda}\cos\left(\frac{1}{\sqrt{\lambda}}\right) - c\lambda \pi
    & =  \frac{2a}{\pi} \int_0^1 \sin(\pi s) \sin\left(\frac{s}{\sqrt{\lambda}}\right)\dd s
         + \frac{2c}{\pi} \int_0^1 (\sin(\pi s))^2 \dd s \\
    & = \frac{a}{\pi}
         \left(\frac{1}{\pi-\frac{1}{\sqrt{\lambda}}} \sin\left(\pi-\frac{1}{\sqrt{\lambda}}\right)
               - \frac{1}{\pi+\frac{1}{\sqrt{\lambda}}} \sin\left(\pi+\frac{1}{\sqrt{\lambda}}\right)
          \right)
        +\frac{c}{\pi},
 \end{align*}
 and hence
 \begin{align}\label{help14}
  \begin{split}
  & \left( \sqrt{\lambda}
          \cos\left(\frac{1}{\sqrt{\lambda}}\right)
          - \frac{1}{\pi(\pi-\frac{1}{\sqrt{\lambda}})} \sin\left(\pi-\frac{1}{\sqrt{\lambda}}\right)
          + \frac{1}{\pi(\pi+\frac{1}{\sqrt{\lambda}})} \sin\left(\pi+\frac{1}{\sqrt{\lambda}}\right)
    \right) a \\
  & -\left(\lambda \pi + \frac{1}{\pi}\right)c
    = 0, \qquad \lambda>0, \quad \lambda\ne \frac{1}{\pi^2}.
  \end{split}
 \end{align}
Note that the eigenfunction \ $e$ \ is identically zero if and only if \ $a=c=0$.
\ Indeed, the functions \ $\sin\left(\frac{t}{\sqrt{\lambda}}\right)$, $t\in[0,1]$,
 \ and \ $\sin(\pi t)$, $t\in[0,1]$, \ are linearly independent provided that \ $\lambda\ne\frac{1}{\pi^2}$, \
  as we check below.
On the contrary let us suppose that they are linearly dependent, i.e., there exist
 constants \ $a,c\in\RR$ \ such that \ $a^2+c^2\ne 0$ \ and \ $a\sin\left(\frac{t}{\sqrt{\lambda}}\right) + c\sin(\pi t) =0$,
 $t\in[0,1]$.
\ Without loss of generality, one can assume that \ $a\ne0$.
\ Then \ $\sin\left(\frac{t}{\sqrt{\lambda}}\right) = -\frac{c}{a}\sin(\pi t)$, $t\in[0,1]$, \ and, by
 differentiating twice, we have \ $\sin\left(\frac{t}{\sqrt{\lambda}}\right) = -\frac{c}{a}\lambda \pi^2 \sin(\pi t)$,
 $t\in[0,1]$.
\ Hence \ $\frac{c}{a}(1-\lambda\pi^2)\sin(\pi t)$, $t\in[0,1]$.
\ Since \ $\lambda\ne\frac{1}{\pi^2}$, \ we have \ $c=0$ \ yielding
 \ $\sin\left(\frac{t}{\sqrt{\lambda}}\right) = 0$, $t\in[0,1]$, \ which is a contradiction.
The system \eqref{help6} and \eqref{help14} has a nontrivial solution \ $(a,c)\ne(0,0)$ \
 if and only if its determinant is zero, which yields
 \begin{align*}
 \pi\lambda^{3/2} \cos\left(\frac{1}{\sqrt{\lambda}}\right)
      + \frac{1}{\pi^2(\pi-\frac{1}{\sqrt{\lambda}})}
         \sin\left(\pi - \frac{1}{\sqrt{\lambda}}\right)
      - \frac{1}{\pi^2(\pi+\frac{1}{\sqrt{\lambda}})}
         \sin\left(\pi + \frac{1}{\sqrt{\lambda}}\right)
      = 0, \qquad \lambda>0, \quad \lambda\ne \frac{1}{\pi^2},
  \end{align*}
 which is equivalent to \eqref{eq_eigenvalue}, since
 \begin{align*}
   \frac{1}{\pi^2(\pi-\frac{1}{\sqrt{\lambda}})}
         \sin\left(\pi - \frac{1}{\sqrt{\lambda}}\right)
      - \frac{1}{\pi^2(\pi+\frac{1}{\sqrt{\lambda}})}
         \sin\left(\pi + \frac{1}{\sqrt{\lambda}}\right)
  & = \frac{1}{\pi^2}\left( \frac{\sin\left(\frac{1}{\sqrt{\lambda}}\right)}{\pi-\frac{1}{\sqrt{\lambda}}}
                          + \frac{\sin\left(\frac{1}{\sqrt{\lambda}}\right)}{\pi+\frac{1}{\sqrt{\lambda}}}
                    \right) \\
  & = \frac{2}{\pi\left(\pi^2 - \frac{1}{\lambda}\right)} \sin\left(\frac{1}{\sqrt{\lambda}}\right).
 \end{align*}

Finally, we check \eqref{eq_eigenfunction}.
By \eqref{help5}, using \ $b=d=0$ \ and \ $\lambda\ne\frac{1}{\pi^2}$, \ we have
 \begin{align*}
  (1-\lambda\pi^2)c
   = 2\int_0^1\sin(\pi s)e(s)\,\dd s
   =a\left( \frac{\sin\left(\pi-\frac{1}{\sqrt{\lambda}}\right)}{\pi-\frac{1}{\sqrt{\lambda}}}
            -  \frac{\sin\left(\pi+\frac{1}{\sqrt{\lambda}}\right)}{\pi+\frac{1}{\sqrt{\lambda}}}  \right)
     + c,
 \end{align*}
 and consequently
 \begin{align*}
  c  = \frac{a}{\lambda \pi^2}
         \left( \frac{\sin\left(\pi+\frac{1}{\sqrt{\lambda}}\right)}{\pi+\frac{1}{\sqrt{\lambda}}}
              -  \frac{\sin\left(\pi-\frac{1}{\sqrt{\lambda}}\right)}{\pi-\frac{1}{\sqrt{\lambda}}}  \right)
      = - \frac{2a}{\lambda \pi \left(\pi^2 - \frac{1}{\lambda}\right)} \sin\left(\frac{1}{\sqrt{\lambda}}\right),
 \end{align*}
 which yields
 \begin{align*}
   e(t) = a\left[ \sin\left(\frac{t}{\sqrt{\lambda}}\right)
                   -\frac{2}{\lambda\pi\left(\pi^2 - \frac{1}{\lambda}\right)}  \sin\left(\frac{1}{\sqrt{\lambda}}\right) \sin(\pi t) \right],
                    \qquad  t\in[0,1].
 \end{align*}
Using \ $-\frac{2}{\lambda\pi\left(\pi^2 - \frac{1}{\lambda}\right)} \sin\left(\frac{1}{\sqrt{\lambda}}\right)
 = \sqrt{\lambda}\pi\cos\left(\frac{1}{\sqrt{\lambda}}\right)$, \ one can finish the proof as in the first
 proof of Corollary \ref{KLcor1}.
\proofend

\bigskip

\noindent{\bf Proof for Remark \ref{Rem_roots}.}
The left-hand side of \eqref{eq_eigenvalue}, denoted by \ $f(\lambda)$, \ is a continuous function
 of \ $\lambda>0$ \ with the extension \ $0$ \ at \ $\lambda=\frac{1}{\pi^2}$ \ due to the fact that
 \begin{align*}
  \lim_{\lambda\to\frac{1}{\pi^2}} f(\lambda)
    = \frac{1}{\pi^3}\cos(\pi) + \frac{2}{\pi^2} \lim_{\lambda\to\frac{1}{\pi^2}}
      \left( -\frac{\sqrt{\lambda}}{2} \cos\left(\frac{1}{\sqrt{\lambda}}\right)\right)
    = -\frac{1}{\pi^3} + \frac{1}{\pi^3}
    =0,
 \end{align*}
 where the first equality follows by L'Hospital's rule.
Let us define a new function \ $\widetilde f\colon(0,\infty)\to\RR$ \ by
 \[
  \widetilde f(\mu):=f\left(\frac{1}{\mu^2}\right)
        =\left\{\begin{array}{ll}
            	\dfrac{1}{\mu^3}\cos(\mu)
                 -\dfrac{2}{\pi^2\left(\mu^2-\pi^2\right)}\sin(\mu)
                        &\mbox{if }\ \mu>0,\,\mu\ne\pi,\\
		            	0&\mbox{if }\ \mu=\pi.
                \end{array}\right.
\]
Since \ $(0,\infty)\ni\mu\mapsto \frac{1}{\mu^2}$ \ is a decreasing homeomorphism of the interval \ $(0,\infty)$ \ onto itself,
  to verify that the equation \eqref{eq_eigenvalue} has a unique root in every interval
  \ $\left(\frac{4}{(2k+1)^2\pi^2}, \frac{4}{(2k-1)^2\pi^2} \right)$, $k\geq 2$, \ $k\in\NN$,
 \ it suffices to prove that \ $\widetilde f$ \ has a unique root in every interval \ $\left(k\pi-\frac \pi2,k\pi+\frac \pi2\right)$,
 $k\geq 2$, \ $k\in\NN$.
\ If \ $\mu\in\left(k\pi-\frac \pi2,k\pi\right)$, \ then \ $\sign(\cos(\mu))=(-1)^k$, \ $\sign(\sin(\mu))=(-1)^{k+1}$,
 \ $\sign(\mu^2 - \pi^2) = 1$, \ thus \ $\sign(\widetilde f(\mu))=(-1)^k$ \ and \ $\widetilde f$ \ has no root in
  \ $\left(k\pi-\frac \pi2,k\pi\right)$, $k\geq 2$, $k\in\NN$.
\ Further, \ $\sign\left(\widetilde f\left(k\pi\right)\right)=\sign(\cos(k\pi))=(-1)^k$ \ and
 \ $\sign\left(\widetilde f\left(k\pi+\frac\pi 2\right)\right)=-\sign\left(\sin\left(k\pi+\frac\pi 2\right)\right)=(-1)^{k+1}$,
 \ thus \ $\widetilde f$ \ has at least one root in \ $\left(k\pi,k\pi+\frac \pi2\right)$, $k\geq 2$, $k\in\NN$.
\ It remains to show that it does not have more than one root.
To this end, we calculate its derivative:
 \begin{align*}
  \widetilde f'(\mu)
   &=-\frac{3}{\mu^4}\cos(\mu)-\frac{1}{\mu^3}\sin(\mu)+\frac{4\mu}{\pi^2(\mu^2-\pi^2)^2}\sin(\mu)-\frac{2}{\pi^2(\mu^2-\pi^2)}\cos(\mu)\\
   &=-\left(\frac{3}{\mu^4}+\frac{2}{\pi^2(\mu^2-\pi^2)}\right)\cos(\mu)
     +\frac{(-\pi^2+4)\mu^4+2\pi^4\mu^2-\pi^6}{\pi^2\mu^3(\mu^2-\pi^2)^2}\sin(\mu)
 \end{align*}
 for \ $\mu\ne \pi$.
\ If \ $\mu\in\left(k\pi,k\pi+\frac \pi2\right)$, $k\geq 2$, $k\in\NN$, \ then the coefficients of \ $\cos(\mu)$ \
 and \ $\sin(\mu)$ \ in the above expression for \ $\widetilde f'(\mu)$ \ are negative.
In case of the coefficient of \ $\sin(\mu)$ \ it follows from \ $\pi^2\mu^3(\mu^2-\pi^2)^2 > 0$ \ and
 \ $(-\pi^2+4)\mu^4+2\pi^4\mu^2-\pi^6 <0$ \ for \ $\mu>2\pi$.
\ Using that \ $\sign(\cos(\mu)) = \sign(\sin(\mu)) = (-1)^k$ \ if \ $\mu\in\left(k\pi,k\pi+\frac \pi2\right)$, $k\geq 2$, $k\in\NN$,
 \ we get \ $\sign(\widetilde f'(\mu)) = (-1)^{k+1}$ \ if \ $\mu\in\left(k\pi,k\pi+\frac \pi2\right)$, $k\geq 2$, $k\in\NN$,
 \ thus \ $\widetilde f$ \ is either strictly increasing or strictly decreasing on this interval,
 and consequently it can have at most one root in this interval.
Finally, the  function \ $f$ \ has no zero greater than \ $\frac{4}{\pi^2}$, \ since
 \begin{align*}
   \lambda^{3/2} \cos\left(\frac{1}{\sqrt{\lambda}}\right)
     + \frac{2}{\pi^2(\pi^2 - \frac{1}{\lambda})}\sin\left(\frac{1}{\sqrt{\lambda}}\right)
    > \frac{2^3}{\pi^3} \cos\left(\frac{\pi}{2}\right)=0,
    \qquad \lambda\in\left(\frac{4}{\pi^2},\infty\right).
 \end{align*}
\proofend

\bigskip

\noindent{\bf Proof of Corollary \ref{KLcor2}.}
One can apply Theorem \ref{THM_KL_gen} with the function \ $g:[0,1]\to\RR$, \ $g(t):=t$, \ $t\in[0,1]$.
\ In the proof of Theorem \ref{THM_KL_gen} we have already checked that the equation
 \ $\sin\left(\frac{1}{\sqrt{\lambda}}\right)=0$ \ is nothing else but the equation \eqref{KL_gen_eigenvalue}
 in the special case \ $g(t) = t$, $t\in[0,1]$.
\ Further, taking into account \ $\lambda=\frac{1}{(k\pi)^2}$, \ $k\in\NN$, \ by partial integration,
 we obtain that the normed eigenfunctions \eqref{KL_gen_eigenvector} take the form
 \begin{align*}
   e(t) &= C\Bigg[ \sqrt{\lambda} \cos\left(\frac{1}{\sqrt{\lambda}}\right) t
                  + \left( \sqrt{\lambda} \sin\left(\frac{1}{\sqrt{\lambda}}\right)
                           +\lambda \left( \cos\left(\frac{1}{\sqrt{\lambda}}\right) - 1 \right)
                    \right) \cos\left(\frac{1}{\sqrt{\lambda}}\right)
                           \sin\left(\frac{t}{\sqrt{\lambda}}\right) \\
       &\phantom{= C\Big[\;}
                   + \left( - \sqrt{\lambda}\cos\left(\frac{1}{\sqrt{\lambda}}\right)
                   + \lambda \sin\left(\frac{1}{\sqrt{\lambda}}\right)
                    \right) \sin\left(\frac{1}{\sqrt{\lambda}}\right)
                           \sin\left(\frac{t}{\sqrt{\lambda}}\right)\\
       &\phantom{= C\Big[\;}
                  + \cos\left(\frac{1}{\sqrt{\lambda}}\right)
                    \cos\left(\frac{t}{\sqrt{\lambda}}\right)
                    \int_0^t u \sin\left(\frac{u}{\sqrt{\lambda}}\right)\dd u \\
       &\phantom{= C\Big[\;}
                  - \cos\left(\frac{1}{\sqrt{\lambda}}\right)
                    \sin\left(\frac{t}{\sqrt{\lambda}}\right)
                    \int_0^t u \cos\left(\frac{u}{\sqrt{\lambda}}\right)\dd u
          \Bigg] \\
       &=  C(-1)^k\Bigg[ \frac{1}{k\pi}t
                        + \frac{(-1)^k - 1}{(k\pi)^2}\sin(k\pi t)
                        + \cos(k\pi t) \int_0^t u \sin(k\pi u)\,\dd u
                        - \sin(k\pi t) \int_0^t u \cos(k\pi u)\,\dd u
                   \Bigg]\\
       &=\frac{C}{(k\pi)^2}\sin(k\pi t),    \qquad t\in[0,1],
 \end{align*}
 where \ $C\in\RR$ \ is such that \ $\int_0^1 (e(t))^2\,\dd t=1$.
\ Since \ $\int_0^1 \sin^2(k\pi t)\,\dd t=\frac{1}{2}$, \ $k\in\NN$, \ we have \ $C=\pm\sqrt{2}(k\pi)^2$, \
 yielding \ $e(t) = \pm\sqrt{2}\sin(k\pi t)$, \ $t\in[0,1]$, \ i.e., we have \eqref{eq_eigenfunction_WB}.
\proofend

\appendix

\vspace*{5mm}

\noindent{\bf\Large Appendix}

\section{ Connections with the paper \cite{Naz2} of Nazarov}\label{Appendix_Nazarov}

In this appendix we compare the Gauss process given in \eqref{def_g_wiener_bridge}
 with the Gauss process given in (1.3) in Nazarov \cite{Naz2}, and then we also compare
 our Theorem \ref{THM_KL_gen} with the results in Section 3 in Nazarov \cite{Naz2}
 for the KL expansions of the Gauss processes in question.

Let \ $(X_t)_{t\in[0,1]}$ \ be a zero-mean Gauss process with continuous sample paths almost surely,
 and suppose that its covariance function \ $G_X$ \ is continuous on \ $[0,1]\times[0,1]$.
\ Let \ $\varphi:[0,1]\to\RR$ \ be a measurable function such that
 \ $\int_0^1 \vert \varphi(s)\vert\,\dd s<\infty$.
\ Introduce the function \ $\psi:[0,1]\to\RR$,
 \begin{align}\label{Nazarov_psi}
  \psi(t):=\int_0^1 G_X(t,s) \varphi(s)\,\dd s, \qquad t\in[0,1].
 \end{align}
Then, by the dominated convergence theorem, \ $\psi$ \ is continuous.
Indeed, if \ $t_n\to t$ \ as \ $n\to\infty$, \ where \ $t_n\in[0,1]$, \ $n\in\NN$, \ and \ $t\in[0,1]$,
 \ then, by the continuity of \ $G_X$, \ for all \ $s\in[0,1]$, \ we have
 \ $\lim_{n\to\infty} G_X(t_n,s) = G_X(t,s)$, \ and there exists \ $K_X>0$ \ such that
 \ $\vert G_X(u,v)\vert\leq K_X$ \ for all \ $(u,v)\in[0,1]\times[0,1]$.
\ Let us denote
 \begin{align}\label{Nazarov_q}
    q:=\int_0^1 \psi(t)\varphi(t)\,\dd t
      = \int_0^1 \int_0^1 G_X(t,s) \varphi(t)\varphi(s) \,\dd t\,\dd s <\infty.
 \end{align}
The finiteness of \ $q$ \ follows from the fact that \ $\psi$ \ (or \ $G_X$) \ is continuous.
For all \ $\alpha\in\RR$, \ introduce the stochastic process
 \begin{align}\label{Nazarov_process}
  X^{\varphi,\alpha}_t := X_t - \alpha \psi(t) \int_0^1 X_s \varphi(s)\,\dd s, \qquad t\in[0,1],
 \end{align}
 where the Lebesgue integral is well-defined almost surely due to the fact that \ $X$ \ has continuous sample paths almost surely.
Then \ $(X^{\varphi,\alpha}_t)_{t\in[0,1]}$ \ is a zero-mean Gauss process with continuous sample paths almost surely and
 with covariance function
 \begin{align}\label{Nazarov_covariance}
  \begin{split}
   G_{X^{\varphi,\alpha}}(t,s)
   & = \cov\left( X_t - \alpha \psi(t) \int_0^1 X_u \varphi(u)\,\dd u,
                 X_s - \alpha \psi(s) \int_0^1 X_v \varphi(v)\,\dd v \right) \\
   & = \cov(X_t,X_s) - \alpha \psi(s) \int_0^1 \cov(X_t,X_u)\varphi(u)\,\dd u  \\
   &\phantom{=\;}
     - \alpha \psi(t) \int_0^1 \cov(X_v,X_s)\varphi(v)\,\dd v \\
   &\phantom{=\;}
      + \alpha^2 \psi(t) \psi(s) \int_0^1 \int_0^1 \cov(X_u,X_v) \varphi(u)\varphi(v) \,\dd u\,\dd v  \\
   &= G_X(t,s) + (q\alpha^2 - 2\alpha)\psi(t)\psi(s),
        \qquad t,s\in[0,1],
 \end{split}
 \end{align}
 see also Nazarov \cite[formula (1.4)]{Naz2}.

Nazarov \cite[Section 3]{Naz2} gave a procedure for finding the KL expansion of \ $(X^{\varphi,\alpha}_t)_{t\in[0,1]}$ \
 supposing that we know the KL expansion of \ $(X_t)_{t\in[0,1]}$.
\ He also provided several examples (specifying \ $X$ \ and \ $\varphi$), \ where he made the KL expansion of
  \ $(X^{\varphi,\alpha}_t)_{t\in[0,1]}$ \ as explicit as possible.
The process \ $X^{\varphi,\alpha}$ \ can be considered as a one-dimensional linear perturbation of the Gauss process \ $X$.

Now we consider two special cases of the above construction.

Let \ $(X_t)_{t\in[0,1]}$ \ be a standard Wiener process, \ $\varphi:[0,1]\to\RR$ \ be a continuous function such that
 \[
   q  = \int_0^1 \int_0^1 (t\wedge s) \varphi(t)\varphi(s)\,\dd t\,\dd s\in(0,1]
   \qquad  \text{and} \qquad  \int_0^1 \left(\int_t^1\varphi(s)\,\dd s\right)^2\dd t = 1,
 \]
 and let \ $\alpha:=\frac{1\pm \sqrt{1-q}}{q}$.
\ For example, one can choose \ $\varphi(t)=\frac{1}{\sqrt{2}\pi} \cos(\pi t)$, $t\in[0,1]$.
\ Then \ $G_X(t,s) = t\wedge s$, $t,s\in[0,1]$, \ $\psi:[0,1]\to\RR$,
 \[
  \psi(t)= \int_0^1 (t\wedge s)\varphi(s)\,\dd s
         = \int_0^t s\varphi(s)\,\dd s + t \int_t^1 \varphi(s)\,\dd s,
         \qquad t\in[0,1],
 \]
 and \ $q\alpha^2 - 2\alpha=-1$.
\ Further, \ $\psi'(t) = \int_t^1 \varphi(s)\,\dd s$, $t\in[0,1]$, \ and \ $\psi''(t) = -\varphi(t)$, $t\in[0,1]$,
 \ satisfying \ $\psi(0) = 0$, \ $\psi'(1)=0$, \ and \ $\int_0^1 (\psi'(s))^2\,\dd s=1$.
\ Hence, taking into account Proposition \ref{Pro1} and \eqref{Nazarov_covariance},
 choosing \ $g$ \ as the function \ $\psi$, \ the Gauss process \ $(Y_t)_{t\in[0,1]}$
 \ given in \eqref{def_g_wiener_bridge} coincides in law with the Gauss process
 \ $(X^{\varphi,\alpha}_t)_{t\in[0,1]}$ \ given in \eqref{Nazarov_process}.

Let \ $(X_t)_{t\in[0,1]}$ \ be a standard Wiener process, and  \ $g:[0,1]\to\RR$ \ be a
 twice continuously differentiable function with \ $g(0) = 0$, \ $g'(1)=0$, \ and \ $\int_0^1 (g'(u))^2\,\dd u=1$.
Let us define the Gauss process  \ $(X^{\varphi,\alpha}_t)_{t\in[0,1]}$ \ given in \eqref{Nazarov_process}
 with \ $\varphi:= -g''$ \ and \ $\alpha:=1$.
\ Since
 \[
   \int_0^1 (t\wedge s) (-g''(s))\,\dd s = g(t), \qquad t\in[0,1],
 \]
 and
 \[
   \int_0^1\int_0^1 (t\wedge s) (-g''(t))(-g''(s))\,\dd t\,\dd s = 1, \qquad t\in[0,1],
 \]
 by \eqref{Nazarov_psi} and \eqref{Nazarov_q}, we have \ $\psi = g$, \ $q=1$ \ and \ $q\alpha^2 - 2\alpha=-1$.
\ Hence the Gauss process \ $(X^{\varphi,\alpha}_t)_{t\in[0,1]}$ \ given in \eqref{Nazarov_process}
 coincides in law with the Gauss process  \ $(Y_t)_{t\in[0,1]}$ \ given in \eqref{def_g_wiener_bridge}.

Based on the above discussion, if \ $g:[0,1]\to\RR$ \ is a twice continuously differentiable function with
 \ $g(0) = 0$, \ $g'(1)=0$, \ and \ $\int_0^1 (g'(u))^2\,\dd u=1$, \ then the Gauss process \ $(Y_t)_{t\in[0,1]}$ \
 given in \eqref{def_g_wiener_bridge} coincides in law with one of the Gauss processes introduced in Nazarov \cite[formula (1.3)]{Naz2}.

\section*{Acknowledgement}
\noindent We would like to thank Endre Igl\'oi for providing us with the formal motivation of the process
 defined in \eqref{def_g_wiener_bridge} presented in the Introduction.
We are grateful to Yakov Nikitin for drawing our attention to the paper \cite{Naz2} of Nazarov,
 and sending us some recent papers on explicit Karhunen--Lo\`{e}ve expansions as well.
Special gratitude to Alexander Nazarov for explaining us in detail how the setup in his paper \cite{Naz2}
 can be specialized to our case, and for pointing out a mistake in an earlier version of the paper related to the fact that
 zero can be an eigenvalue of the integral operator under consideration.

\bibliographystyle{plain}

\end{document}